%% file: main.tex
\documentclass[pdflatex,sn-mathphys-num]{sn-jnl}

\input{preamble}

\theoremstyle{thmstyleone}%
\theoremstyle{thmstyletwo}%

\theoremstyle{thmstylethree}%

\newcommand\algo{\textsc{2DD-AG2m}\xspace}
\newcommand\sadam{\textsc{AG2m}\xspace}
\newcommand\sada{\textsc{AG2}\xspace}
\newcommand\ddada{\textsc{DD-AG2m}\xspace}

\begin{document}

\title[Two-level AdaGrad for training GNNs]{Two-level domain-decomposition AdaGrad method for scalable training of graph neural networks}


\author[1,2]{\fnm{Laurynas} \sur{Varnas}}\email{laurynas.varnas@toulouse-inp.fr}
\author[1,6]{\fnm{Julien} \sur{Herrmann}}\email{julien.herrmann@irit.fr}
\author[4]{\fnm{Alexander} \sur{Heinlein}}\email{A.Heinlein@tudelft.nl}
\author[1,2,3]{\fnm{Serge} \sur{Gratton}}\email{serge.gratton@toulouse-inp.fr}
\author*[1,2,3,5]{\fnm{Alena} \sur{Kopani\v{c}\'akov\'a}}\email{Alena.Kopanicakova@unige.ch}

\affil[1]{\orgname{Institut de Recherche en Informatique de Toulouse}, \orgaddress{\city{Toulouse}, \country{France}}}
\affil[2]{\orgname{Artificial and Natural Intelligence Toulouse Institute}, \orgaddress{\city{Toulouse}, \country{France}}}
\affil[3]{\orgname{Toulouse-INP (ENSEEIHT)}, \orgaddress{\city{Toulouse}, \country{France}}}
\affil[4]{\orgname{Delft Institute of Applied Mathematics, Delft University of Technology}, \orgaddress{\city{Delft}, \country{Netherlands}}}
\affil[5]{\orgname{Department of Computer Science, University of Geneva}, \orgaddress{\city{Geneva}, \country{Switzerland}}}
\affil[6]{\orgname{Centre National de la Recherche Scientifique}, \orgaddress{\country{France}}}


\abstract{
	Graph neural networks (GNNs) have emerged as a powerful framework for learning from graph-structured data.
	However, their efficient training remains challenging, particularly in distributed computing environments.
	This challenge arises from the use of message passing, which couples all graph nodes, leading to expensive optimization steps, high memory requirements, and substantial communication overhead.
	To alleviate these limitations, we propose a novel domain-decomposition (DD) variant of \sadam, an AdaGrad method enhanced with second-order curvature information and momentum, denoted by \ddada.
	The proposed \ddada alternates between \sadam optimization on the original (global) graph and \sadam optimization on the partitioned graphs.
	To incorporate global information at reduced cost, we further introduce a two-level variant (\algo) that performs global optimization steps on a coarse graph obtained by randomly subsampling nodes within each subdomain.
	Numerical experiments spanning graph classification, node-level regression, and spatiotemporal forecasting tasks demonstrate that the proposed DD methods reduce the computational cost required to achieve the same predictive performance by a factor of $4$--$8$.
	Moreover, for the fixed computational cost, they improve the predictive performance of GNNs by up to $22\%$ compared with the baseline \sadam.
}

\keywords{GNN, training, AdaGrad, domain-decomposition, multilevel optimization}



\maketitle

\input{sections/01_intro.tex}
\input{sections/02_GNN.tex}
\input{sections/03_OFFO.tex}
\input{sections/04_num_examples.tex}
\input{sections/05_num_results.tex}
\input{sections/06_conclusion.tex}

\bmhead{Acknowledgements}
Authors thank A.~Kssim for many insightful discussions and for providing the initial implementation of \sadam algorithm.
This work benefited from ANITI (Grant No.~ANR-23-IACL-0002) and the DAIMOS project (Grant No.~ANR-25-EXNU-0002).
The numerical results were carried out using HPC resources from GENCI-IDRIS (Grant No.~AD011015766R1).

\bibliography{biblio}

\end{document}

%% file: preamble.tex
\usepackage[T1]{fontenc}
\usepackage[utf8]{inputenc}
\usepackage{caption}
\usepackage{algorithm,algpseudocode}
\usepackage{graphicx,url}
\usepackage{amsmath,amssymb,amsopn}
\usepackage{tabularx}
\usepackage{color}
\usepackage{calc}
\usepackage{caption}
\usepackage{longtable}
\usepackage{multirow}
\usepackage{mathrsfs}
\usepackage{array}
\usepackage{hyperref}
\usepackage{cleveref}
\usepackage{tikz}
\usepackage{booktabs}
\usepackage{pgfplotstable}
\usepackage{pgfplots}
\usepackage{amsfonts}
\usepackage{amsthm}
\usepackage[title]{appendix}
\usepackage{textcomp}
\usepackage{manyfoot}
\usepackage{physics}
\usepackage{xfrac}
\usepackage{bm}
\usetikzlibrary{backgrounds,positioning,shapes,calc,arrows.meta}
\usepackage[outline]{contour}
\usepackage{calrsfs}
\usepackage{graphicx}
\usepgfplotslibrary{groupplots}
\usetikzlibrary{matrix}
\usetikzlibrary{automata, external, shapes.geometric, fit}
\def \wv{\vec{w}}

\def \yv{\vec{y}}

\def \rm{\mat{r}}

\def \R{\mathbb{R}}			

\DeclareMathAlphabet{\pazocal}{OMS}{zplm}{m}{n}

\renewcommand{\vec}[1]{\boldsymbol{#1}}
\newcommand{\mat}[1]{\boldsymbol{{#1}}}

\algnewcommand\algorithmiconput{\textbf{Constants:}}
\algnewcommand\algorithmicinput{\textbf{Input:}}
\algnewcommand\algorithmicoutput{\textbf{Output:}}
\algnewcommand{\algorithmicgoto}{\textbf{go to}}%

\algnewcommand\Constants{\item[\algorithmiconput]}
\algnewcommand\Input{\item[\algorithmicinput]}%
\algnewcommand\Output{\item[\algorithmicoutput]}%
\algnewcommand{\Goto}[1]{\algorithmicgoto~\ref{#1}}%

\definecolor{myblack}{RGB}{53, 53, 53}
\definecolor{myblue}{RGB}{40, 75, 99}
\definecolor{myred}{RGB}{192, 50, 33}
\definecolor{myyellow}{RGB}{255, 166, 48}
\definecolor{mywhite}{RGB}{240, 237, 238}
\definecolor{mygreen}{RGB}{0, 102, 0}

\definecolor{green1}{RGB}{9, 82, 86}
\definecolor{green2}{RGB}{8, 127, 140}
\definecolor{green3}{RGB}{6, 167, 125}
\definecolor{green4}{RGB}{79, 109, 122}
\definecolor{green5}{RGB}{192, 214, 223}
\definecolor{violet}{RGB}{26,69,131}

\definecolor{checkgreen}{rgb}{0,0.6,0}
\definecolor{phase1}{rgb}{0.008,0.655,1.000}
\definecolor{phase2}{rgb}{0.016,0.75,0.700}
\definecolor{phase3}{rgb}{0.929,0.35,0.700}
\definecolor{icsyellow}{cmyk}{0.00,0.11,0.53,0.00}

\definecolor{blackmy}{RGB}{38, 70, 83}
\definecolor{bluemy}{RGB}{39, 125, 161}
\definecolor{greenmy}{RGB}{42, 167, 143}
\definecolor{yellowmy}{RGB}{233, 196, 106}
\definecolor{brownmy}{RGB}{244, 162, 97}
\definecolor{redmy}{RGB}{249, 65, 68}

\definecolor{darkbluemy}{RGB}{65, 59, 147}
\definecolor{lightbluemy}{RGB}{71, 139, 194}
\definecolor{greenmy}{RGB}{98, 173, 153}
\definecolor{darkorangemy}{RGB}{230, 142, 52}
\definecolor{lightorangemy}{RGB}{217, 172, 59}

\definecolor{blue1}{RGB}{1, 42, 74}
\definecolor{blue2}{RGB}{1, 73, 124}
\definecolor{blue3}{RGB}{42, 111, 151}
\definecolor{blue4}{RGB}{44, 125, 160}
\definecolor{blue5}{RGB}{70, 143, 175}
\definecolor{blue6}{RGB}{137, 194, 217}

\definecolor{red1}{RGB}{204, 68, 75}
\definecolor{red2}{RGB}{218, 85, 82}
\definecolor{red3}{RGB}{227, 150, 149}
\definecolor{red4}{RGB}{228, 190, 171}

\definecolor{brown1}{RGB}{92,178,112}
\definecolor{brown2}{RGB}{130,194,110}
\definecolor{brown3}{RGB}{163, 193, 173}

\definecolor{steelblue}{RGB}{70,130,180}
\definecolor{sandybrown}{RGB}{244,164,96}
\definecolor{crimson}{RGB}{220,20,60}
\definecolor{mediumpurple}{RGB}{147,112,219}
\definecolor{forestgreen}{RGB}{34,139,34}

\usepackage[normalem]{ulem}
\usepackage{soul,xcolor}
\usepackage[many]{tcolorbox}
\usepackage{mathtools}
\usepackage{algorithm}
\usepackage{enumitem}
\usepackage{listings}
\usepackage{relsize}
\usepackage{nameref}
\usepackage{xspace}

\usepackage{longtable}

\newcommand{\algorithmiccommentMine}[1]{\bgroup\hfill$\triangleright$~{ \textcolor{gray}{#1}}\egroup}
\newcommand\COMMENTmine[1]{\algorithmiccommentMine{#1}}

\newcommand\oldtext[1]{}
\newcommand\cancel[1]{}

\newcommand\oldtextt[1]{}

\usetikzlibrary{spy}
\usepackage{pgfplots}
\usepgfplotslibrary{fillbetween}

\usepackage{array}
\newcommand{\PreserveBackslash}[1]{\let\temp=\\#1\let\\=\temp}
\newcolumntype{C}[1]{>{\PreserveBackslash\centering}p{#1}}
\newcolumntype{R}[1]{>{\PreserveBackslash\raggedleft}p{#1}}
\newcolumntype{L}[1]{>{\PreserveBackslash\raggedright}p{#1}}

\DeclareMathAlphabet\bpazocal{OMS}{cmsy}{b}{n}
\usetikzlibrary{patterns}

%% file: sections/01_intro.tex
\section{Introduction}

Graph Neural Networks (GNNs) have emerged as a powerful class of neural
architectures for learning from graph-structured data by combining entity
attributes with relational information encoded by the graph topology~\cite{Scarselli2009GraphNN,Battaglia2018RelationalIB,wu2020comprehensive}.
Following the early recurrent formulation of
Scarselli et al.~\cite{Scarselli2009GraphNN}, several architectural
principles have been introduced, including spectral graph
convolutions~\cite{Defferrard2016ConvolutionalNN,Kipf2017SemiSupervisedCG},
inductive neighborhood aggregation~\cite{Hamilton2018InductiveRL},
attention-weighted aggregation~\cite{Velickovic2018GraphAN}, and higher-order
or substructure-aware models designed to improve expressivity
~\cite{Morris2019WeisfeilerLG,Bouritsas2023ImprovingGN}.
Many of these architectures can be expressed within the message-passing (MP)
framework~\cite{Gilmer2017NeuralMP}, in which the graph topology plays a dual
role: it describes the relational structure of the data and determines the
computational dependencies through which information is propagated.

GNNs have been applied to a wide range of problems in which the relations
between entities carry information that cannot be captured from their
attributes alone, such as molecular
property prediction and atomistic modeling~\cite{Gilmer2017NeuralMP,Schutt2021EquivariantMP}, social and recommender systems~\cite{Hamilton2018InductiveRL}, and knowledge-graph reasoning~\cite{Schlichtkrull2018ModelingRD,Chami2020LowDimensionalHK}.
Of particular interest are scientific applications in which the graph represents the discretization of a physical or spatiotemporal system. For instance, GNNs have been used as graph-based surrogate models for physical simulation~\cite{SanchezGonzalez2020LearningSC}, mesh-based numerical simulation~\cite{Pfaff2021LearningMS}, computational fluid dynamics~\cite{Bonnet2022AirfRANSHF}, and traffic forecasting~\cite{Li2018DiffusionCR}.
Moreover, weather forecasting provides a prominent large-scale spatiotemporal application, with GraphCast~\cite{lam2022graphcast} and GenCast~\cite{price2025gencast} employing mesh-based GNNs for deterministic and probabilistic global forecasting, respectively.
In such applications, increasing the spatial resolution produces larger and denser graphs,
and therefore raises the computational and memory costs
of both MP and backpropagation.

Despite the success of GNNs accros wide range of applications, training GNNs efficiently at scale remains a significant challenge.
In contrast to dense feed-forward networks, MP layers combine sparse neighborhood aggregation with dense
transformations of node or edge representations.
Irregular graph topologies and
non-uniform node degrees therefore induce
indirect memory accesses and workload imbalance, which can limit accelerator utilization and complicate parallel
execution~\cite{Jia2020RedundancyFreeCG,wang2021gnnadvisor}.
When training is distributed across multiple devices, MP dependencies additionally require exchanging node features and gradients between workers.
The irregularity and volume of these exchanges hinder communication–computation overlap, increase synchronization costs, and may ultimately dominate the training time~\cite{wan2022bns,Shao2024DistributedGN}.

These memory and computational demands have motivated a variety of methods that operate on reduced computation graphs.
Node- and layer-wise sampling restrict the neighborhoods explored to a set of target nodes~\cite{Hamilton2018InductiveRL,Zou2019LayerDependentIS}, whereas subgraph-based methods construct mini-batches from connected graph regions, either through graph partitioning, as in Cluster-GCN~\cite{Chiang2019ClusterGCNEA}, or subgraph sampling with normalized stochastic estimators, as in GraphSAINT~\cite{Zeng2019GraphSAINTGS}.
Distributed GNN systems additionally address graph placement, remote feature access, workload balancing, and inter-device communication~\cite{Shao2024DistributedGN,wan2022pipegcn,wan2022bns}.
These approaches primarily reorganize or approximate the forward and backward computations.

Closer to our setting are methods based on local training over graph partitions.
For example, LLCG~\cite{ramezani2021learn} trains independent copies of a GNN on graph partitions, applies periodic parameter averaging, and uses global server corrections to compensate for the dependencies omitted during local training.
Its analysis shows that periodic averaging alone may
retain a residual error when cross-partition dependencies are ignored.
Randomized partitioning has alternatively been proposed to improve model
aggregation by reducing the discrepancy between the data distributions,
losses, and gradients observed by the local
trainers~\cite{zhu2025simplifying}. Other approaches preserve part of the
missing graph information directly during MP; for instance,
DIGEST~\cite{chai2022distributed} periodically synchronizes stale
representations of remote neighbors.

Meanwhile, the numerical analysis community has extensively developed domain-decomposition (DD) methods~\cite{Quarteroni1999DomainDM, Toselli2005DomainDM, Saad2003IterativeMS} for solving, in parallel, large-scale systems of equations arising from the discretization of partial differential equations (PDEs).
The idea behind these methods is to decompose the original problem into smaller subproblems that can be solved concurrently, while coordinating the information exchange among them in order to recover a globally coherent solution.
Owing to this design, DD methods provide a natural framework for parallelization, and they rank among the most scalable solution strategies in scientific computing available to date.

The interplay between DD and machine learning has recently attracted significant attention; see~\cite{Klawonn2024MachineLD,Heinlein2020CombiningML} for a comprehensive overview.
In this work, we focus on approaches for accelerating and parallelizing neural network training.
In this spirit, several methods based on decomposing either the network parameters or the training data have been proposed and led to model-parallel~\cite{Kopanicakova2024EnhancingTP, lee2026two} or data-parallel training~\cite{CruzAlegria2025DataParallel}, respectively.
Moreover, authors  of \cite{Gratton2023MultilevelOFFO, Gratton2025RecursiveBA} introduced a recursive bound-constrained AdaGrad framework, which extends AdaGrad-type methods to the multilevel and DD settings while accounting for stochastic gradient estimates.
Complementarily,~\cite{Gal2025EfficientTG} proposed a multiscale approach to GNN training based on a hierarchy of coarsened graphs.

Motivated by these developments, we propose two DD-based variants of AdaGrad for the parallel training of GNNs.
The proposed methods extend the recursive bound-constrained AdaGrad framework of~\cite{Gratton2025RecursiveBA} to the training of GNNs, with the decomposition induced directly by the structure of the input graph.
The central idea is to partition the input graph into subgraphs and to perform independent local optimization steps on the resulting subgraphs using temporary copies of a shared global model.
The local corrections are subsequently aggregated to update the parameters of the global model, in the spirit of additive Schwarz methods~\cite{Toselli2005DomainDM}.
To incorporate global information at a reduced cost, we further introduce a two-level variant, which performs global optimization steps using a coarse graph, obtained by per-subdomain random node subsampling.
In contrast to conventional distributed GNN training methods, which primarily focus on computational scalability, the proposed approach treats graph partitioning as a nonlinear preconditioning mechanism for the AdaGrad optimizer, thereby enhancing convergence speed while also exposing additional parallelism.

The remainder of the paper is organized as follows.
In Section~\ref{sec:gnn}, we review GNNs and formulate the associated training problem.
Section~\ref{sec:dd_adagrad} introduces the proposed \ddada and \algo algorithms, which constitute the main contribution of this work.
In Section~\ref{sec:num_examples}, we describe the benchmark problems, together with the algorithmic and implementation details.
Section~\ref{sec:num_results} demonstrates the convergence and the algorithmic scalability properties of the proposed algorithms.
Finally, in Section~\ref{sec:conclusion}, we summarize our findings and discuss possible future extensions.

%% file: sections/02_GNN.tex
\section{Graph Neural Networks (GNN)}
\label{sec:gnn}
GNNs are neural architectures designed to process data whose underlying
structure is represented by a graph~\cite{wu2020comprehensive,Battaglia2018RelationalIB}.

\subsection{The Graph Neural Network Model}
The input data of a GNN are represented by an attributed graph.
Formally, a graph is defined as a pair $\pazocal{G} = (\pazocal{V}, \pazocal{E})$, where $\pazocal{V}$ denotes the set of nodes and $\pazocal{E} \subset \pazocal{V} \times \pazocal{V}$ denotes the set of edges connecting pairs of nodes. For any node $v \in \pazocal{V}$, the neighborhood $\pazocal{N}(v)$ is defined as the set of nodes sharing an edge with $v$, i.e.,
$$
	\pazocal{N}(v) = \left\{ u \in \pazocal{V} : (u,v) \in \pazocal{E} \text{ or } (v,u) \in \pazocal{E} \right\}.
$$
The connectivity of $\pazocal{G}$ is encoded by the (possibly weighted) adjacency matrix ${\bm{A} \in \mathbb{R}^{|\pazocal{V}| \times |\pazocal{V}|}}$, whose entry $\bm{A}_{uv}$ is nonzero if and only if $u$ and $v$ are connected.
Each node $v$ carries a feature vector $\bm{x}_v \in \mathbb{R}^{d}$, encoding node-specific attributes (e.g., physical quantities, sensor measurements, or pixel intensities, depending on the application), with $d$ being the number of features per node.
These vectors are collected row-wise into the node feature matrix $\bm{X} \in \mathbb{R}^{|\pazocal{V}| \times d}$.
When available, each edge $(u,v) \in \pazocal{E}$ also carries a feature vector $\bm{e}_{uv}$.

A GNN is a function $f_{\bm{\theta}}(\cdot)$, parameterized by the parameters $\bm{\theta} \in \mathbb{R}^{n}$, which maps an input graph $\pazocal{G}$ to a prediction $f_{\bm{\theta}}(\pazocal{G})$.
The output may be associated with the complete graph, individual nodes, or node-time pairs. For instance, these settings correspond to graph classification, node regression, and spatiotemporal forecasting, respectively, all of which are considered in the numerical experiments of this work, see Section~\ref{sec:num_results}.

\subsection{Message-passing (MP) neural networks}
Many widely used local GNN architectures can be formulated within the
message-passing (MP) framework~\cite{Gilmer2017NeuralMP, Scarselli2009GraphNN},
which provides a flexible way to model dependencies between nodes in
graph-structured data. In this framework, each node updates its representation
by aggregating information from its local neighborhood, with the graph
structure defining the corresponding computational dependencies.

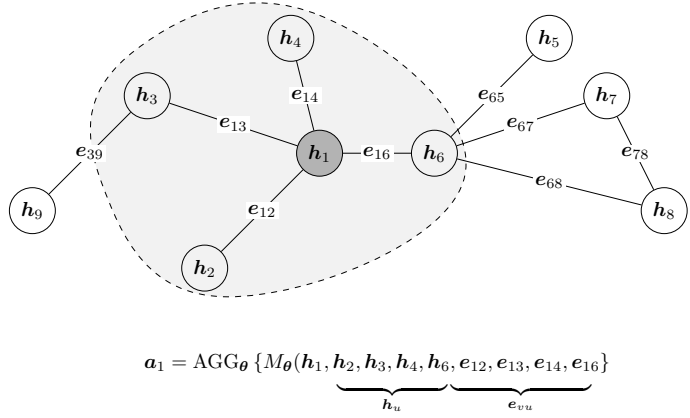
\begin{figure}[tb!]
	\centering
	\resizebox{0.7\textwidth}{!}{\input{fig/MessagePassing}}
	\caption{Illustration of MP aggregation for node $v_1$ with feature vector $\bm{h}_1$. The sketch is adapted from~\cite{Scarselli2009GraphNN}.}
	\label{fig:neighborhood}
\end{figure}

The initial hidden features of node $v$ are given as $\bm h_v^{(0)} = \bm x_v$.
At each layer $\ell=0,\ldots,L-1$,  node representations are updated as
\begin{equation*}
\begin{aligned}
    \bm a_v^{(\ell)}
    &=
    \operatorname{AGG}_{\bm{\theta}, u\in\pazocal N(v)}
    \left\{
        M_{\bm{\theta}}^{(\ell)}
        \left(
            \bm h_v^{(\ell)},
            \bm h_u^{(\ell)},
            \bm e_{uv}
        \right)
    \right\},
    \\
    \bm h_v^{(\ell+1)}
    &=
    U_{\bm{\theta}}^{(\ell)}
    \left(
        \bm h_v^{(\ell)},
        \bm a_v^{(\ell)}
    \right),
\end{aligned}
\end{equation*}
where $\bm{h}_v^{(\ell)}$ is the hidden state of node $v$ at layer $\ell$.
The functions $M_{\bm{\theta}}^{(\ell)}$ and $U_{\bm{\theta}}^{(\ell)}$, parametrized by $\bm{\theta}$, are the learnable message and differentiable update functions, respectively.
Moreover, $\operatorname{AGG}_{\bm{\theta}}$ denotes a permutation-invariant aggregation operator, such as element-wise sum, mean, maximum, or attention-weighted sum.
Following the MP formulation of~\cite{Gilmer2017NeuralMP}, the message function determines the information
passed from a neighboring node to $v$ based on the representations of the two nodes and, when present, the features of their connecting edge.
The aggregation then combines these messages to form the input to the node update.
Since each layer propagates information by at most one hop, the representation
\(\bm{h}_v^{(L)}\) obtained after \(L\) layers depends only on the information contained in the
\(L\)-hop neighborhood of \(v\), as depicted in Figure~\ref{fig:neighborhood}.

The node representations obtained after $L$ layers of MP, given by $\bm{h}_v^{(L)}$, are mapped to task-specific predictions through a readout function $R_{\bm{\theta}}$.
This defines the overall parameterized map
$$
	f_{\bm{\theta}}(\pazocal{G}) = R_{\bm{\theta}}\qty(\qty{\bm{h}_v^{(L)} : v \in \pazocal{V}}).
$$
The specific form of $R_{\bm{\theta}}$ depends on the task at hand.
For node-level tasks, $R_{\bm{\theta}}$ maps each final representation
$\bm{h}^{(L)}_v$ independently to a node-wise prediction.
For graph-level tasks, the representations
$\{\bm{h}^{(L)}_v : v \in \pazocal{V}\}$ are first combined through a permutation-invariant pooling operator, such as sum or mean, before being mapped to a single graph-level prediction.
Spatiotemporal architectures may additionally include recurrent or temporal-processing layers, but their graph component still follows the local MP structure described above.

Note that, in the above description of MP, the vector $\bm{\theta} \in \mathbb{R}^{n}$ denotes the collection of all learnable parameters of the message, update, aggregation, and readout functions.
The parameters of the message and update functions are shared across all nodes in the graph.
Because the aggregation operator is invariant to the ordering of the incoming messages and the same functions are applied at every node, the resulting node representations are equivariant to permutations of the nodes.
Furthermore, the MP framework imposes no restriction on the size of a node’s neighborhood.

The MP formulation includes many commonly used architectures, for instance
Graph Convolutional Networks~\cite{Kipf2017SemiSupervisedCG}, GraphSAGE~\cite{Hamilton2018InductiveRL}, and Graph Attention Networks~\cite{Velickovic2018GraphAN} mainly differ in their definitions of
$M_{\bm{\theta}}^{(\ell)}$, $U_{\bm{\theta}}^{(\ell)}$, and $\operatorname{AGG}_{\bm{\theta}}$.

\subsection{GNN training}
We consider a dataset $\pazocal{D} = \qty{(\pazocal{G}^{(j)}, \bm{y}^{(j)})}_{j=1}^{N}$ of $N$ samples,
where $\pazocal{G}^{(j)}$ is an attributed input graph, and $\bm{y}^{(j)}$ contains its supervised target.
Given a task-specific loss function that compares the prediction $f_{\bm{\theta}}(\pazocal{G}^{(j)})$ with the corresponding target $\bm{y}^{(j)}$, training consists of minimizing the following empirical risk:
\begin{equation}
	\min_{\bm{\theta} \in \R^n} \pazocal{L}(\bm{\theta}) := \frac{1}{N} \sum_{j=1}^{N} {\mathrm{loss}}\qty(f_{\bm{\theta}}(\pazocal{G}^{(j)}), \bm{y}^{(j)}).
	\label{eq:empirical_risk}
\end{equation}

The cost of solving~\eqref{eq:empirical_risk} depends on the number of samples, as well as the size and the topology of the input graphs.
Each MP layer combines sparse neighborhood aggregation with dense transformations of node or edge representations.
Non-uniform node degrees consequently induce irregular workloads and memory-access patterns, which may reduce accelerator utilization and complicate parallel execution~\cite{Jia2020RedundancyFreeCG,wang2021gnnadvisor}.

These costs are further increased by gradient-based optimization, since backpropagation requires storing or recomputing intermediate node representations and MP dependencies across GNN layers.
Moreover, increasing the number of MP layers enlarges the receptive field of each node at the cost of increasing the amount of intermediate information involved in both the forward and backward computations.
The resulting memory footprint therefore scales as $O(L\,|\pazocal{V}|\,d')$, growing with the network depth $L$, the size of the active graph $|\pazocal{V}|$, and the width $d'$ of the hidden representations~\cite{li2021training}.

Furthermore, when a training sample is a very large graph, the memory required for full-graph processing may exceed the available device memory.
Sampling, clustering, and graph partitioning address this limitation by restricting computations to smaller subgraphs~\cite{Chiang2019ClusterGCNEA,Zeng2019GraphSAINTGS,Hamilton2018InductiveRL}.
This enables independent and potentially parallel processing of different subgraphs.
However, restricting MP to subgraphs changes the loss function and its gradient.
In contrast, the DD-based algorithms proposed in this work use subgraph and subsampled computations solely to construct local corrections, which are subsequently coordinated through global optimization steps to ensure consistency with the original full-graph training problem.

%% file: fig/MessagePassing.tex
\tikzsetnextfilename{message_passing}
\begin{tikzpicture}[
		node/.style={circle, draw, minimum size=8mm, inner sep=0pt},
		shadednode/.style={circle, draw, fill=black!30, minimum size=8mm, inner sep=0pt},
		edgelabel/.style={midway, fill=white, inner sep=1pt}
	]

	\node[shadednode] (n1) at (0,0)    {$\bm{h}_1$};
	\node[node]       (n2) at (-2,-2)  {$\bm{h}_2$};
	\node[node]       (n3) at (-3,1)   {$\bm{h}_3$};
	\node[node]       (n4) at (-0.5,2) {$\bm{h}_4$};
	\node[node]       (n6) at (2,0)    {$\bm{h}_6$};
	\node[node]       (n5) at (4,2)    {$\bm{h}_5$};
	\node[node]       (n7) at (5,1)    {$\bm{h}_7$};
	\node[node]       (n8) at (6,-1)   {$\bm{h}_8$};
	\node[node]       (n9) at (-5,-1)  {$\bm{h}_9$};

	\begin{scope}[on background layer]
		\fill[black!10, opacity=0.5]
		plot[smooth cycle, tension=0.8] coordinates {
				(-3.2,1.8)
				(0,2.4)
				(2.5,-0.5)
				(-1.6,-2.5)
				(-3.8,-0.5)
			};
		\draw[dashed]
		plot[smooth cycle, tension=0.8] coordinates {
				(-3.2,1.8)
				(0,2.4)
				(2.5,-0.5)
				(-1.6,-2.5)
				(-3.8,-0.5)
			};
	\end{scope}

	\draw (n1) -- node[edgelabel] {$\bm{e}_{12}$} (n2);
	\draw (n1) -- node[edgelabel] {$\bm{e}_{14}$} (n4);
	\draw (n1) -- node[edgelabel] {$\bm{e}_{13}$} (n3);
	\draw (n1) -- node[edgelabel] {$\bm{e}_{16}$} (n6);
	\draw (n3) -- node[edgelabel] {$\bm{e}_{39}$} (n9);
	\draw (n6) -- node[edgelabel] {$\bm{e}_{65}$} (n5);
	\draw (n6) -- node[edgelabel] {$\bm{e}_{67}$} (n7);
	\draw (n6) -- node[edgelabel] {$\bm{e}_{68}$} (n8);
	\draw (n7) -- node[edgelabel] {$\bm{e}_{78}$} (n8);

	\node[align=center] at (1,-4)
	{$\bm{a}_1 = \operatorname{AGG}_{\bm{\theta}} \qty{M_{\bm{\theta}}(\bm{h}_1, \bm{h}_2, \bm{h}_3, \bm{h}_4, \bm{h}_6, \bm{e}_{12}, \bm{e}_{13}, \bm{e}_{14}, \bm{e}_{16}}$ \\[-0.7em]
		$\phantom{a_1 = \operatorname{AGG}_{\bm{\theta}} M_{\bm{\theta}}(h_1,} \underbrace{\phantom{h_2, h_3, h_4, h_6}}_{\bm{h}_u} \underbrace{\phantom{\bm{e}_{13}, \bm{e}_{14}, \bm{e}_{13}, \bm{e}_{16} } }_{\bm{e}_{vu}}$};

\end{tikzpicture}

%% file: sections/03_OFFO.tex
\section{Two-level DD-based AdaGrad for GNNs}
\label{sec:dd_adagrad}

In this work, we minimize~\eqref{eq:empirical_risk} using AdaGrad like algorithms from the family of \emph{objective-function-free optimization} (OFFO) methods~\cite{Gratton2025NoiseTolerant, Gratton2023MultilevelOFFO, Gratton2025RecursiveBA}.
OFFO methods rely solely on (approximate) gradient information and, when available, (approximate) second-order information, without requiring objective function evaluations.
This property makes them robust to noise introduced by subsampling.

\subsection{Single-level curvature-informed AdaGrad with momentum}
We begin our presentation of the proposed algorithms by first discussing the single-level variant of AdaGrad, enhanced with curvature information~\cite{Gratton2025NoiseTolerant} and momentum, which we abbreviate as \sadam.
This algorithm constitutes the core computational engine of the proposed \ddada and \algo methods.
In contrast to the AdaGrad variant considered in~\cite{Gratton2025NoiseTolerant}, we also incorporate a momentum term.

\begin{algorithm}[tb!]
	\small
	\caption{\sadam($K$, $\pazocal{L}$, $\bm{\theta}_0$, $\bm{w}_0$, $\bm{m}_0$, $adjust$)}
	\label{alg:adagrad2}
	\begin{algorithmic}[1]
		\Statex \textbf{Input:} $K \in \mathbb{N}$; $\pazocal{L}:\R^n \to \R $; $\bm{\theta}_0 \in \R^{n}$;  $\bm{w}_0 \in \R^{n}_{\ge 0}$; $\bm{m}_0 \in \R^{n}$; $adjust \in \left\{ true, false \right\}$
		\Statex \textbf{Constants:} $\beta \in [0,1)$

		\For{$k = 1$ to $K$}

		\State $\bm{g}_{k} \approx \grad_{\bm{\theta}} \pazocal{L} (\bm{\theta}_{k-1})$ \COMMENTmine{Get (stochastic) gradient estimate}

		\State $\bm{w}_{k, i} = \sqrt{\bm{w}_{k-1, i}^{\,2} + \bm{g}_{k, i}^{\,2}}$ \COMMENTmine{Compute AdaGrad weights}

		\If{$k=1$ and $adjust$}
		\State $\bm{w}_{k, i} = \max\qty(\bm{w}_{k, i}, \bm{w}_{0, i})$ \COMMENTmine{Adjust weights based on the full-space} \label{line:weight_adjust}
		\EndIf

		\State $\bm{\Delta}_{k, i} = \dfrac{ |\bm{g}_{k, i}| }{ \bm{w}_{k, i} }$ \COMMENTmine{Evaluate TR-radius}

		\State $\bm{s}^{S}_{k, i} = \min\qty( \max\qty( -\bm{g}_{k, i}, -\bm{\Delta}_{k, i}), \bm{\Delta}_{k, i})$ \COMMENTmine{Compute CP step}

		\State $\bm{B}_{k} \approx \grad_{\bm{\theta}}^2 \pazocal{L} (\bm{\theta}_{k-1})$ \COMMENTmine{Get Hessian approximation}

		\State
		$\bm{s}^{Q}_{k} = \gamma_{k} \bm{s}^{S}_{k}$,  where \COMMENTmine{Get second-order informed step-size}
		\begin{align}
			\gamma_{k} =
			\begin{cases}
				\min\!\qty(
				1,
				\dfrac{
					-\langle\bm{g}_{k},
					\bm{s}^{S}_{k}\rangle
				}{
					\langle \bm{s}^{S}_{k},
					\bm{B}_{k}
					\bm{s}^{S}_{k} \rangle
				}
				),
				   & \text{if } \langle \bm{s}^{S}_{k},  \bm{B}_{k} \bm{s}^{S}_{k} \rangle >0, \\
				1, & \text{otherwise}
			\end{cases}
			\label{eq:QP_sub}
		\end{align}

		\State $\bm{m}_{k} = \beta \bm{m}_{k-1} + (1-\beta)\bm{s}^{Q}_{k}$ \COMMENTmine{Evaluate momentum}

		\State $\bm{m}_{k, i} = \min\qty( \max\qty( \bm{m}_{k, i}, -\bm{\Delta}_{k, i}), \bm{\Delta}_{k, i})$ \COMMENTmine{Enforce TR bound}

		\State $\bm{\theta}_{k} = \bm{\theta}_{k-1} + \bm{m}_{k}$ \COMMENTmine{Update parameters}

		\EndFor

		\State \Return $\bm{\theta}_{K}, \bm{w}_{K}, \bm{m}_{K}$
	\end{algorithmic}
\end{algorithm}

As outlined in Algorithm~\ref{alg:adagrad2}, at the $k$-th iteration, \sadam first evaluates the gradient or its subsampled approximation, i.e., $\bm{g}_{k} \approx \nabla \pazocal{L}(\bm{\theta}_{k-1})$.
The accumulated squared gradients are then used to compute the coordinate-wise AdaGrad weights $\bm{w}_{k}$.
These weights are then used to define the trust-region (TR) radius
$\bm{\Delta}_{k,i} = |\bm{g}_{k,i}| / \bm{w}_{k,i}$, where the subscript $i$ denotes the $i$-th component of a vector.
Following the terminology of TR methods~\cite{Conn2000TrustRegion}, the vector $\bm{\Delta}_{k}$ defines a coordinate-wise, adaptive bound on the step size.
Thus, if the accumulated weights increase, the corresponding admissible step sizes decrease, in turn reducing the magnitude of the admissible updates along coordinates that repeatedly exhibit large gradients.
The Boolean flag $adjust$ indicates whether the initial AdaGrad weights must be adjusted using the previously accumulated global (full-space) weights.
As we will see later, this adjustment is performed whenever \sadam is invoked on a coarse or subgraph problem, ensuring that the size of their first \sadam step is bounded by the global weights.

Once $\bm{\Delta}_{k}$ has been computed, \sadam first determines the step $\bm{s}^{S}_{k}$ obtained by truncating the negative gradient direction to the TR.
The quality of this search direction is then improved by scaling it with a step-size $\gamma_{k} \in \R$, computed from (approximate) second-order information.
Specifically, $\gamma_{k}$ is computed such that $\gamma_{k}\bm{s}^{S}_{k}$ minimizes the local quadratic model $\langle \bm{g}_{k}, \bm{s}\rangle + \tfrac12 \langle \bm{s}, \bm{B}_{k}\bm{s}\rangle$ along the direction $\bm{s}^{S}_{k}$, where $\bm{B}_{k}$ is either Hessian $\grad^2 \pazocal{L}(\bm{\theta}_{k-1})$ or its symmetric approximation.
Note that if curvature is negative or the resulting step size is larger than one, the direction $\bm{s}^{S}_{k}$ is kept.

Finally, the step $\bm{s}^{Q}_{k} = \gamma_{k}\bm{s}^{S}_{k}$ is incorporated into the momentum as ${\bm{m}_{k} = \beta\bm{m}_{k-1} + (1-\beta)\bm{s}^{Q}_{k}}$.
This momentum is subsequently bounded to ensure that the resulting update satisfies
$-\bm{\Delta}_{k,i} \le \bm{m}_{k,i} \le \bm{\Delta}_{k,i}$.
The vector $\bm{m}_{k}$ is then used to update the GNN parameters, yielding the new iterate $\bm{\theta}_{k}$.
The adopted momentum scheme follows the same algorithmic structure as in~\cite{erway2020trust, kopanicakova2022globally}.
Moreover, it falls within the class of methods analyzed recently in~\cite{Gratton2026UnifiedAdaptive}, where the global convergence rate is established for a broad family of adaptive first-order methods with momentum.

\subsection{DD-based AdaGrad for GNNs}
Training GNNs on large-scale graphs is computationally demanding, as each parameter update requires MP over the entire graph.
To improve scalability and exploit parallelism, we adopt a DD strategy, well known from solving large-scale systems associated with discretized PDEs~\cite{Toselli2005DomainDM}.
In particular, we propose \ddada, a DD-based variant of the \sadam~\cite{Gratton2023MultilevelOFFO, Gratton2025RecursiveBA} method, specifically designed for training GNNs.
The proposed \ddada method partitions the input graph into multiple subgraphs, enabling optimization steps to be performed independently and in parallel on each subgraph.
The resulting subgraph updates are subsequently synchronized to produce an update for the GNN parameters.

Following the DD terminology, throughout the description of the \ddada algorithm we refer to computations on subgraphs as \emph{subdomain computations} and computations on the original graph as \emph{global computations}.
Superscripts are used to distinguish quantities associated with the global graph or subgraphs.
Moreover, we use double subscripts to denote the outer and inner iterations of the \ddada algorithm.
For instance, $\bm{\theta}^p_{k,K}$ denotes the model parameters associated with the $p$-th subdomain after the $k$-th outer (\ddada) iteration and the $K$-th inner (\sadam) iteration.

\subsubsection{Global and subdomain minimization problems}
\label{sec:problems_def}

The \ddada is associated with the minimization of global and subgraph loss functions.
The global loss function $\pazocal{L}^G$ is defined using original graphs, i.e.,~$\pazocal{L}^G:=\pazocal{L}$.

\begin{figure}[tb!]
	\centering
	\resizebox{0.8\textwidth}{!}{
		\input{fig/GraphPartition}}
	\caption{An example of non-overlapping graph partitioning. The original
		graph (left) and a partition of the node set $\pazocal{V}$ into $P = 3$
		disjoint subsets $\pazocal{V}^{1}, \pazocal{V}^{2}, \pazocal{V}^{3}$
		(right).}
	\label{fig:graph_partition}
\end{figure}
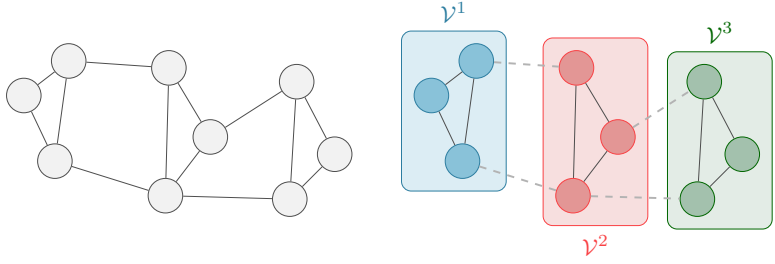

The subdomain losses are defined on the partitioned graphs.
Thus, given a graph $\pazocal{G}$, we partition the node set $\pazocal{V}$ into $P$ disjoint subsets $\{\pazocal{V}^1,\ldots,\pazocal{V}^P\}$ such that
\[
	\pazocal{V}^i \cap \pazocal{V}^j = \emptyset \quad (i \neq j),
	\qquad
	\bigcup_{p=1}^P \pazocal{V}^p = \pazocal{V}.
\]
Each partition induces a subgraph $\pazocal{G}^p = (\pazocal{V}^p, \pazocal{E}^p)$, where $\pazocal{E}^p$ contains only edges whose endpoints both belong to $\pazocal{V}^p$, i.e., all inter-partition edges are removed\footnote{Although our formulation assumes non-overlapping partitions, the proposed framework can be naturally extended to overlapping graph partitions.}.
An illustration of the non-overlapping graph partitioning is provided in Figure~\ref{fig:graph_partition}.

Each subgraph is associated with a subdomain loss, given as
\begin{equation}
	\pazocal{L}^{p}(\bm{\theta}^{p})
	:= \frac{1}{N} \sum_{j=1}^N
	{\mathrm{loss}}\qty(f_{\bm{\theta}^{p}}(\pazocal{G}^{(p,j)}), \yv^{(p,j)}),
	\qquad p = 1,\ldots,P,
	\label{eq:local_losses}
\end{equation}
where the upper script $(p,j)$ denotes $p$-th partition of $j$-th sample.
All subdomain losses $\{\pazocal{L}^{p}\}_{p=1}^{P}$ are defined using the same GNN architecture as the global problem~\eqref{eq:empirical_risk}.
Consequently, each subdomain maintains its own parameter vector $\bm{\theta}^{p}\in\mathbb{R}^{n}$, which has the same structure as the global parameter vector $\bm{\theta}$, since only the graph is decomposed, while the model parameters are not.

\subsubsection{\ddada algorithm}
Having defined subdomain losses~\eqref{eq:local_losses}, we now describe the \ddada method, outlined in Algorithm~\ref{alg:dd}.
Starting from the current iterate $\bm{\theta}_{k-1,0}$, each iteration of \ddada consists of three stages.
First, \sadam is applied to the global loss $\pazocal{L}^{G}:=\pazocal{L}$ for $K^G$ inner steps, producing the iterate $\bm{\theta}_{k-1,K^G}$.
This global step operates on the entire graph $\pazocal{G}$ and therefore ensures global information propagation.

Afterward, the DD step is performed on the partitioned graphs by applying
\sadam to minimize each subdomain loss $\pazocal{L}^{p}$ for $K^p$ inner steps.
The subdomain optimization process is initialized using a copy of the global parameters, obtained after minimization of the global loss, i.e., $\bm{\theta}_{k-1, K^G}$.
Moreover, the AdaGrad weights are initialized by utilizing the weights $\wv_{k,K^G}$, computed after the global optimization step.
In particular, the initial subdomain weights are constructed by taking the maximum between the weights constructed during the first subdomain step and the global AdaGrad weights (i.e., \sadam is invoked with $adjust = true$).
This ensures that the first subdomain correction satisfies
$|\bm{s}_{1,i}| \leq \bm{\Delta}_{1,i} \leq |\bm{g}_{1,i}|/\bm{w}_{0,i}$, for each component $i$, thereby bounding its magnitude using the corresponding global AdaGrad weights, see~\cite{Gratton2025RecursiveBA}.
Since the subdomain minimization problems are independent of each other, this DD phase can be carried out in parallel.

After the DD phase is completed, the obtained correction is computed as
${\bm{s}^{p}_{k} = \bm{\theta}^{p}_{k,K^G+K^p} - \bm{\theta}_{k-1,K^G}}$, for $p = 1, \ldots, P$.
These corrections are then reconciled into a single global update, denoted by $\bm{s}^{DD}_{k} \in \R^n$.
Here, we construct $\bm{s}^{DD}_{k} \in \R^n$ by means of averaging, i.e., $\bm{s}^{DD}_{k} = \frac{1}{P} \sum_{p=1}^{P} \bm{s}^{p}_{k}$.
Afterwards, the parameters of global GNN, $\bm{\theta}_{k-1,K^G}$, are updated as
${\bm{\theta}_{k,0} = \bm{\theta}_{k-1,K^G} + \bm{s}^{DD}_{k}}$.
This requires only a single communication step, consisting of a global reduction of all corrections across the subdomains.
In the end, we point out that more sophisticated aggregation strategies have recently been proposed in~\cite{salvado2026multi} and could be incorporated into the proposed DD framework in future work.

The algorithm then proceeds to the next outer iteration.

\begin{algorithm}[tb!]
	\small
	\caption{\ddada ($\pazocal{L}^G, \qty{\pazocal{L}^{p}}_{p=1}^{P}, \bm{\theta}_{0,0}, \bm{w}_{0,0}, \bm{m}_{0}$)}
	\label{alg:dd}
	\begin{algorithmic}[1]
		\Statex \textbf{Input:}
		$\pazocal{L}^G: \R^n \to \R$;
		$\qty{\pazocal{L}^{p}}_{p=1}^{P}$, where each $\pazocal{L}^p: \R^n \to \R$;
		$\bm{\theta}_{0,0} \in \R^{n}$,
		$\bm{w}_{0,0} \in \R^{n}_{\ge 0}$,
		$\bm{m}_{0} \in \R^{n}$
		\Statex \textbf{Constants:} $P, K, K^G, \qty{K^p}_{p=1}^{P} \in \mathbb{N}$

		\For{$k = 1$ to $K$}

		\COMMENTmine{Perform $K^G$ \sadam steps using the fine graph}
		\State $\bm{\theta}_{k-1,K^G}, \bm{w}_{k-1,K^G}, \bm{m}_{k-1,K^G} =$
		\Call{\sadam}{$K^G, \pazocal{L}^{G}, \bm{\theta}_{k-1,0}, \bm{w}_{k-1,0}, \bm{m}_{k-1,0}, false$}

		\State
		\For{$p = 1$ to $P$ \textbf{ in parallel}}

		\State
		\COMMENTmine{Perform $K^p$ \sadam steps using $p$-th subgraph}
		\State $\bm{\theta}^{p}_{k,K^G+K^p}, \ldots =$
		\Call{\sadam}{$K^p, \pazocal{L}^{p}, \bm{\theta}_{k-1,K^G}, \bm{w}_{k-1,K^G}, \bm{m}_{k-1,K^G}, true$}

		\State $\bm{s}^{p}_{k} = \bm{\theta}^{p}_{k,K^G+K^p} - \bm{\theta}_{k-1,K^G}$ \COMMENTmine{Evaluate subdomain correction}

		\EndFor

		\State $\bm{\theta}_{k,0} = \bm{\theta}_{k-1,K^G} +  \sum_{p=1}^P \frac{1}{P} \bm{s}^{p}_{k}$ \COMMENTmine{Update global GNN with DD corrections}

		\EndFor

		\State \Return $\bm{\theta}_{K,0}$

	\end{algorithmic}
\end{algorithm}

\subsection{Two-level DD-based AdaGrad for GNNs}
The \ddada algorithm exploits subgraph computations to accelerate \sadam's convergence while enabling parallelism.
Its main computational bottleneck is the $K^G$ \sadam steps performed on the global graph $\pazocal{G}$, which are inherently serial.
To reduce this cost, we replace several global optimization steps with optimization steps on a coarsened graph.
This follows the standard two-level DD approaches~\cite{dolean2015introduction}, in which a coarse level efficiently transfers global information at reduced cost.

\subsubsection{Coarse (subsampled graph) minimization problem}
To define the coarse minimization problem, we first construct the coarse graph.
Following~\cite{Gal2025EfficientTG}, we express the coarsening process through a transfer operator that maps the original graph data onto a coarse representation.
Concretely, let $\bm{R} \in \{0,1\}^{n_{\mathrm{c}} \times |\pazocal{V}|}$ denote a binary restriction (subsampling) operator, whose rows are the indicator vectors of the $n_{\mathrm{c}}$ retained coarse-level nodes.
The retained nodes are selected by random subsampling with coarsening factor $c_f$, performed independently within each subdomain, while ensuring that at least one node is retained from every subdomain.

\begin{figure}[tb!]
	\centering
	\input{fig/GraphCoarsening}
	\caption{Sketch of the per-subdomain random node subsampling to create the coarse graph.
		The original graph (left) and the coarse graph $\pazocal{G}^{C}$, formed by the selected nodes together with the edges inherited between them (right).}
	\label{fig:graph_coarsening}
\end{figure}
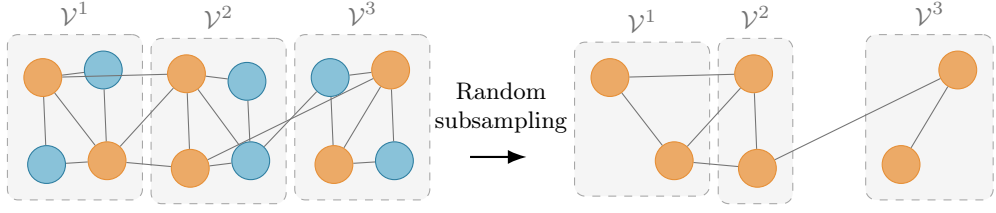

Thus, for every subdomain $\pazocal{V}^{p}$, a random subset of nodes is kept and the remaining nodes are dismissed.
The operator $\bm{R}$ then restricts the global graph node set to the coarse node set as $\bm{X}^{C} = \bm{R}\,\bm{X}$.
This defines the coarse graph $\pazocal{G}^{C} = (\pazocal{V}^{C}, \pazocal{E}^{C})$ with node features $\bm{X}^{C}$.
The edges are thus directly inherited from $\pazocal{G}$.
Using standard Galerkin projection, well-known from multigrid, the adjacency matrix $\bm{A}^{C}$ is constructed as $\bm{A}^{C} = \bm{R}\,\bm{A}\,\bm{R}^{\top}$.
Moreover, for node-level tasks, the labels are restricted analogously, $\bm{y}^{C} = \bm{R}\,\bm{y}$.
Note, performing the random node subsampling per subdomain ensures that the coarsening respects the partitioning. At the same time, the retained inter-subdomain edges keep $\pazocal{G}^{C}$ globally connected and therefore reflective of the global structure of graph $\pazocal{G}$.
An illustration of the employed coarsening strategy is given in Figure~\ref{fig:graph_coarsening}.

The coarse-level loss function $\pazocal{L}^{C}: \R^{n} \to \R$ is defined as
\begin{equation}
	\pazocal{L}^{C}(\bm{\theta}^C)
	:= \frac{1}{N} \sum_{j=1}^N
	{\mathrm{loss}} (f_{\bm{\theta}^{C}}(\pazocal{G}^{(C,j)}), \yv^{(C,j)}),
	\label{eq:coarse_loss}
\end{equation}
where the superscript $(C,j)$ denotes $j$-th coarse sample.
Similarly to subdomain losses, $\pazocal{L}^{C}$ is evaluated using the GNN with the same architecture as the global problem~\eqref{eq:empirical_risk}.
Thus, the transfer operator $\bm{R}$ acts only on the graph data, i.e., the node features, adjacency matrix, and labels.
The parameters of global GNN, $\bm{\theta}$, are not coarsened, thus $\bm{\theta}^C \in \mathbb{R}^{n}$.
This is possible as the weight matrices act on the feature dimension rather than on the node dimension, so that the same MP operations can be applied on the global and the coarse graphs~\cite{Gal2025EfficientTG}.

We note that~\cite{Gal2025EfficientTG} considers additional coarsening
strategies that could be explored within the proposed DD framework in future
work.

\subsection{\algo algorithm}
Algorithm~\ref{alg:tldd} summarizes the proposed two-level \ddada algorithm, denoted as \algo.
As we can see, \algo follows the same structure as \ddada.
However, the global steps are alternating, not just with subdomain steps, but also with coarse steps.
Here, we point out that the convergence of the \algo is controlled by the full-space AdaGrad weights.
To this end, \sadam is invoked with $adjust = true$ during both the coarse and the subdomain minimization, so that the initial step sizes during coarse and subdomain minimization are restricted by the global AdaGrad weights. 

\begin{algorithm}[tb!]
	\small
	\caption{\algo ($\pazocal{L}^G, \pazocal{L}^C, \qty{\pazocal{L}^{p}}_{p=1}^{P}, \bm{\theta}_{0,0}, \bm{w}_{0,0}, \bm{m}_{0}$)}
	\label{alg:tldd}
	\begin{algorithmic}[1]
		\Statex \textbf{Input:}
		$\pazocal{L}^G := \pazocal{L}: \R^n \to \R$;
		$\pazocal{L}^C:  \R^n \to \R$;
		$\qty{\pazocal{L}^{p}}_{p=1}^{P}$, where each $\pazocal{L}^p: \R^n \to \R$;
		$\bm{\theta}_{0,0} \in \R^{n}$,
		$\bm{w}_{0,0} \in \R^{n}_{\ge 0}$,
		$\bm{m}_{0} \in \R^{n}$
		\Statex \textbf{Constants:} $P, K, K^{G}, K^C, \qty{K^p}_{p=1}^{P} \in \mathbb{N}$

		\For{$k = 1$ to $K$}

		\COMMENTmine{Perform $K^{G}$ \sadam steps using the original graph}
		\State $\bm{\theta}_{k-1,K^{G}}, \bm{w}_{k-1,K^{G}}, \bm{m}_{k-1,K^{G}} =$
		\Call{\sadam}{$K^{G}, \pazocal{L}^{G}, \bm{\theta}_{k-1,0}, \bm{w}_{k-1,0}, \bm{m}_{k-1,0}, false$}

		\COMMENTmine{Perform $K^C$ \sadam steps using the coarse graph}
		\State $\bm{\theta}_{k-1,K^{G} + K^C}, \_ , \_ =$
		\Call{\sadam}{$K^C, \pazocal{L}^{C}, \bm{\theta}_{k-1,K^{G}}, \bm{w}_{k-1,K^{G}}, \bm{m}_{k-1,K^{G}}, true$}

		\Statex

		\COMMENTmine{Perform $K^{G}$ \sadam steps using the original graph}
		\State $\bm{\theta}_{k-1,\,2K^{G}+K^C}, \bm{w}_{k-1,\,2K^{G}+K^C}, \bm{m}_{k-1,\,2K^{G}+K^C} =$
		\Statex \hfill \Call{\sadam}{$K^{G}, \pazocal{L}^{G}, \bm{\theta}_{k-1,K^{G}+K^C}, \bm{w}_{k-1,K^{G}+K^C}, \bm{m}_{k-1,K^{G}}, false$}

		\Statex

		\For{$p = 1$ to $P$ \textbf{ in parallel}}

		\COMMENTmine{Perform $K^p$ subdomain \sadam steps using $p$-th subgraph}
		\State $\bm{\theta}^{p}_{k,\,2K^{G}+K^C+K^p}, \_ , \_ =$
		\Statex \hfill
		\Call{\sadam}{$K^p, \pazocal{L}^{p}, \bm{\theta}_{k-1,\,2K^{G}+K^C}, \bm{w}_{k-1,\,2K^{G}+K^C}, \bm{m}_{k-1,\,2K^{G}+K^C}, true$}

		\State $\bm{s}^{p}_{k} = \bm{\theta}^{p}_{k,\,2K^{G}+K^C+K^p} - \bm{\theta}_{k-1,\,2K^{G}+K^C}$ \COMMENTmine{Evaluate subdomain correction}
		\EndFor

		\State $\bm{\theta}_{k,0} = \bm{\theta}_{k-1,\,2K^{G}+K^C} +  \sum_{p=1}^P \frac{1}{P} \bm{s}^{p}_{k}$ \COMMENTmine{Update global GNN with DD corrections}

		\EndFor

		\State \Return $\bm{\theta}_{K,0}$

	\end{algorithmic}
\end{algorithm}

We also highlight that both \ddada and \algo use momentum during the subdomain and coarse optimization.
Following~\cite{kopanicakova2022globally}, the corresponding momentum variables are initialized from the global momentum, allowing the subdomain/coarse optimizations to benefit from the optimization history associated with the global objective.
However, the subdomain and coarse momentums are not propagated back.
This prevents optimization histories specific to the coarse/subdomain objectives from being mixed with the momentum associated with the global objective.

%% file: fig/GraphPartition.tex
\tikzsetnextfilename{graph_partition}
\begin{tikzpicture}[
		gnode/.style={circle, draw=black!70, fill=black!5, minimum size=5mm, inner sep=0pt},
		intra/.style={black!70},
		inter/.style={gray!60, dashed, thick},
	]


	\begin{scope}
		\foreach \name/\x/\y in {%
				a1/0/2.0, a2/0.65/2.5, a3/0.45/1.05,
				b1/2.1/2.4, b2/2.7/1.4, b3/2.05/0.55,
				c1/3.95/2.2, c2/4.5/1.15, c3/3.85/0.5}
		\node[gnode] (\name) at (\x,\y) {};
		\foreach \u/\v in {%
				a1/a2, a1/a3, a2/a3, b1/b2, b1/b3, b2/b3, c1/c2, c1/c3, c2/c3,
				a2/b1, a3/b3, b2/c1, b3/c3}
		\draw[intra] (\u) -- (\v);
	\end{scope}

	\begin{scope}[xshift=5.9cm]
		\foreach \name/\x/\y in {a1/0/2.0, a2/0.65/2.5, a3/0.45/1.05}
		\node[gnode, draw=blue4, fill=blue6] (\name) at (\x,\y) {};
		\foreach \name/\x/\y in {b1/2.1/2.4, b2/2.7/1.4, b3/2.05/0.55}
		\node[gnode, draw=redmy, fill=red3] (\name) at (\x,\y) {};
		\foreach \name/\x/\y in {c1/3.95/2.2, c2/4.5/1.15, c3/3.85/0.5}
		\node[gnode, draw=mygreen, fill=brown3] (\name) at (\x,\y) {};

		\begin{scope}[on background layer]
			\node[draw=blue4, fill=blue6!30, rounded corners, inner sep=5pt,
				fit=(a1)(a2)(a3), label={[blue4,font=\small]above:$\pazocal{G}^1$}] {};
			\node[draw=redmy, fill=red3!30, rounded corners, inner sep=5pt,
				fit=(b1)(b2)(b3), label={[redmy,font=\small]below:$\pazocal{G}^2$}] {};
			\node[draw=mygreen, fill=brown3!30, rounded corners, inner sep=5pt,
				fit=(c1)(c2)(c3), label={[mygreen,font=\small]above:$\pazocal{G}^3$}] {};
		\end{scope}

		\foreach \u/\v in {%
				a1/a2, a1/a3, a2/a3, b1/b2, b1/b3, b2/b3, c1/c2, c1/c3, c2/c3}
		\draw[intra] (\u) -- (\v);
		\foreach \u/\v in {a2/b1, a3/b3, b2/c1, b3/c3}
		\draw[inter] (\u) -- (\v);

	\end{scope}

\end{tikzpicture}

%% file: fig/GraphCoarsening.tex
\tikzsetnextfilename{graph_coarsening}
\begin{tikzpicture}[
		sel/.style={circle, draw=darkorangemy, fill=darkorangemy!75, minimum size=5mm, inner sep=0pt},
		drop/.style={circle, draw=blue4, fill=blue6, minimum size=5mm, inner sep=0pt},
		hull/.style={draw=black!35, dashed, fill=black!4, rounded corners, inner sep=6pt},
		edge/.style={black!55},
	]

	\begin{scope}
		\node[sel]  (a1) at (0.1,2.3)  {};
		\node[drop] (a2) at (0.9,2.4)  {};
		\node[drop] (a3) at (0.15,1.15){};
		\node[sel]  (a4) at (0.95,1.2) {};
		\node[sel]  (b1) at (2.0,2.35) {};
		\node[drop] (b2) at (2.8,2.25) {};
		\node[sel]  (b3) at (2.05,1.1) {};
		\node[drop] (b4) at (2.85,1.2) {};
		\node[drop] (c1) at (3.9,2.3)  {};
		\node[sel]  (c2) at (4.7,2.4)  {};
		\node[sel]  (c3) at (3.95,1.15){};
		\node[drop] (c4) at (4.75,1.2) {};

		\begin{scope}[on background layer]
			\node[hull, fit=(a1)(a2)(a3)(a4), label={[black!55,font=\small]above:$\pazocal{V}^1$}] {};
			\node[hull, fit=(b1)(b2)(b3)(b4), label={[black!55,font=\small]above:$\pazocal{V}^2$}] {};
			\node[hull, fit=(c1)(c2)(c3)(c4), label={[black!55,font=\small]above:$\pazocal{V}^3$}] {};
		\end{scope}

		\foreach \u/\v in {a1/a2, a1/a3, a2/a4, a3/a4, a1/a4,
				b1/b2, b1/b3, b2/b4, b3/b4, b1/b4,
				c1/c2, c1/c3, c2/c4, c3/c4, c2/c3,
				a1/b1, a4/b1, a4/b3, b3/c2, b4/c1}
		\draw[edge] (\u) -- (\v);

	\end{scope}

	\node[font=\small, align=center] at (6.16,1.9) {Random\\subsampling};
	\draw[-{Latex[length=2.5mm]}, line width=0.8pt] (5.75,1.25) -- (6.5,1.25);

	\begin{scope}[xshift=7.5cm]
		\node[sel] (qa1) at (0.1,2.3)  {};
		\node[sel] (qa4) at (0.95,1.2) {};
		\node[sel] (qb1) at (2.0,2.35) {};
		\node[sel] (qb3) at (2.05,1.1) {};
		\node[sel] (qc2) at (4.7,2.4)  {};
		\node[sel] (qc3) at (3.95,1.15){};

		\begin{scope}[on background layer]
			\node[hull, fit=(qa1)(qa4), label={[black!55,font=\small]above:$\pazocal{V}^1$}] {};
			\node[hull, fit=(qb1)(qb3), label={[black!55,font=\small]above:$\pazocal{V}^2$}] {};
			\node[hull, fit=(qc2)(qc3), label={[black!55,font=\small]above:$\pazocal{V}^3$}] {};
		\end{scope}

		\foreach \u/\v in {qa1/qa4, qb1/qb3, qc2/qc3, qa1/qb1, qa4/qb1, qa4/qb3, qb3/qc2}
		\draw[edge] (\u) -- (\v);

	\end{scope}

\end{tikzpicture}

%% file: sections/04_num_examples.tex
\section{Numerical examples and implementation details}
\label{sec:num_examples}

In this section, we describe the benchmark problems utilized to evaluate the proposed \ddada and \algo algorithms, together with their algorithmic configuration and implementation details.
The selected benchmark problems cover diverse learning tasks, GNN architectures, and graph partitioning strategies, enabling an assessment of the convergence behavior of the proposed DD-based methods across a broad range of settings.
Table~\ref{tab:benchmark_summary} summarizes the characteristics of all benchmark problems and the corresponding algorithmic setup.

\begin{table}[tb!]
	\caption{Summary of the benchmark problems and the associated experimental setup.}
	\label{tab:benchmark_summary}
	\begin{tabular*}{\textwidth}{@{\extracolsep\fill}l|l|l|l}
		\toprule
		& CIFAR10 & AirfRANS & METR-LA \\
		\midrule
		Learning task & Graph classification & Node regression & Spatiotemporal regression \\
		GNN architecture & GCN~\cite{Kipf2017SemiSupervisedCG} & GraphSAGE~\cite{Bonnet2022AirfRANSHF} & DCRNN~\cite{Li2018DiffusionCR} \\
		Number of layers  & $4$ & $4$ & $2$ \\
		Layer width & $146$ & $64$ & $32$ \\
		Loss function & Cross-entropy & MSE & Masked MAE \\
		Graph partitioner & Spectral~\cite{VonLuxburg2007TutorialSC} & Morton ordering~\cite{Morton1966Geodetic} & Multilevel Graclus~\cite{Dhillon2007WeightedGC} \\
		Nodes per graph & $85$--$150$ & $32{,}000$ & $207$\footnotemark[1] \\
		Batch size & $32$ & $1$ & $64$ \\
		\botrule
	\end{tabular*}
\end{table}
\footnotetext[1]{For METR-LA, all samples share a single static graph with $207$ nodes, whereas for CIFAR10 and AirfRANS each sample has its own graph.}

\subsection{Benchmark problems}
We consider the following three benchmark problems:

\subsubsection{Super-pixel CIFAR10}
\begin{figure}[tb!]
	\centering
	\includegraphics[width=0.85\textwidth]{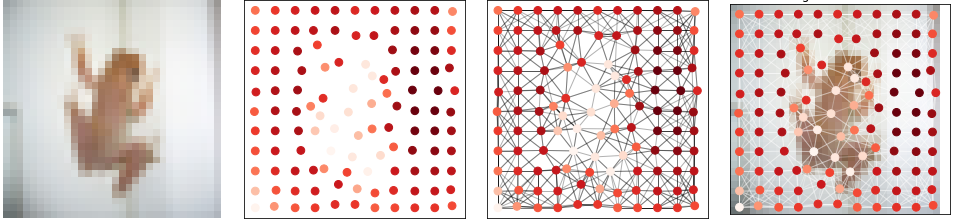}
	\caption{Example of a sample from CIFAR10 and its corresponding super-pixel graph~\cite{Dwivedi2022BenchmarkingGN}: the original image (left), the super-pixel nodes (middle left), the $k$-nearest-neighbor graph (middle right), and the resulting graph (right).}
	\label{fig:cifar_superpixel}
\end{figure}

Our first benchmark problem concerns graph classification on the CIFAR10 super-pixel dataset.
This dataset is derived from the CIFAR10 image dataset~\cite{Krizhevsky2009CIFAR}, which contains $32 \times 32$ RGB images from ten object categories.
In particular, we use the preprocessed super-pixel graphs from the benchmark suite of~\cite{Dwivedi2022BenchmarkingGN}, constructed using the Simple Linear Iterative Clustering (SLIC) algorithm for image segmentation~\cite{knyazev2019understanding}.
SLIC partitions each image into spatially compact regions of pixels with similar color characteristics, which define the graph nodes whose features consist of the mean RGB intensity of the super-pixel and its two-dimensional centroid coordinates.
Edges connect each super-pixel to its eight nearest neighbors and are weighted using a Gaussian kernel, resulting in graphs with approximately $85$--$150$ nodes; see Figure~\ref{fig:cifar_superpixel}.
We use the standard split of $45{,}000$ / $5{,}000$ / $10{,}000$ graphs for training, validation, and testing, respectively.

Graph classification is performed using a graph convolutional network (GCN)~\cite{Kipf2017SemiSupervisedCG}, with the architecture specified in~\cite{Dwivedi2022BenchmarkingGN}, which comprises four graph convolutional layers with a hidden width of $146$, batch normalization, residual connections, and mean readout.
The GCN is trained with the cross-entropy loss.
We partition the graphs using spectral partitioning~\cite{VonLuxburg2007TutorialSC}, for which the adjacency matrix is converted to a dense representation.

\subsubsection{AirfRANS}
\begin{figure}[tb!]
	\centering
	\includegraphics[width=\textwidth]{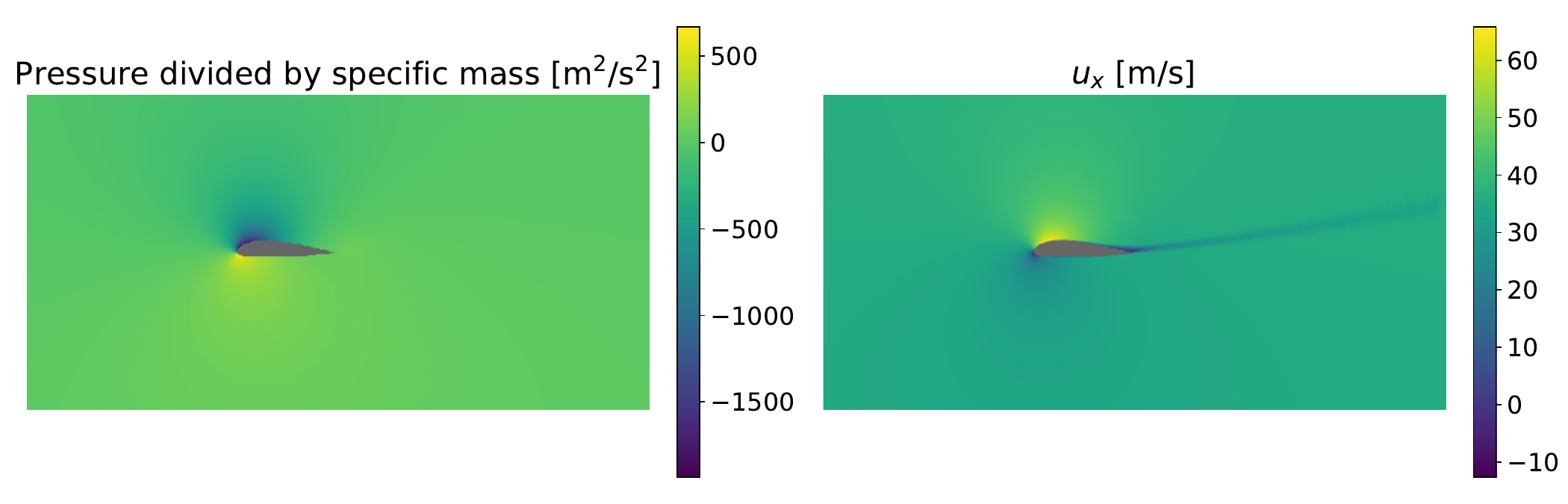}
	\caption{Example of flow fields around a NACA airfoil from the AirfRANS dataset: the pressure divided by specific mass (left) and the velocity field in the $x$-direction (right).}
	\label{fig:airfrans_fields}
\end{figure}

Next, we consider a physics-based regression benchmark in computational fluid dynamics, based on the AirfRANS dataset~\cite{Bonnet2022AirfRANSHF}.
It contains high-fidelity solutions of the two-dimensional incompressible steady-state Reynolds-averaged Navier--Stokes (RANS) equations over airfoils from the NACA 4- and 5-digit families.
Each sample represents the steady turbulent flow around a single airfoil in the subsonic regime for Reynolds numbers in $[2 \times 10^6,\,6 \times 10^6]$ and angles of attack in $[-5^\circ,\,15^\circ]$.
Each node is described by five input features, i.e.,  the two free-stream velocity components, the distance to the airfoil surface, and the two components of the surface normal (set to zero at interior nodes).
The regression targets are the two local velocity components, the reduced pressure (i.e., the pressure divided by the fluid density), and the turbulent kinematic viscosity.
Representative pressure and velocity fields are shown in Figure~\ref{fig:airfrans_fields}.

Following~\cite{Bonnet2022AirfRANSHF}, we consider the scarce-data setting with $200$ training and $200$ testing simulations.
Each simulation is provided as an unstructured point cloud without explicit connectivity.
As in the baseline setup, we uniformly subsample each simulation to $32{,}000$ nodes and construct a sparse geometric graph by connecting nodes whose normalized Euclidean distance is below $0.05$.

We employ the baseline GraphSAGE architecture of~\cite{Bonnet2022AirfRANSHF}, comprising four SAGE convolution layers of width $64$, preceded by an MLP encoder and followed by an MLP decoder.
The network is trained using the mean squared error (MSE) loss.
Graph partitioning is performed using a Morton (Z-order) space-filling curve~\cite{Morton1966Geodetic}, i.e., the normalized point coordinates are converted to Morton codes, sorted, and divided into contiguous segments of equal size.
We choose this partitioning because it relies solely on the spatial coordinates and is therefore independent of the graph connectivity.
To construct the coarse-level dataset, we randomly sample $32{,}000 / c_f$ nodes from the original point cloud and connect them using a radius graph with radius $0.1$.

\subsubsection{METR-LA}
\begin{figure}[tb!]
	\centering
	\includegraphics[width=.45\textwidth]{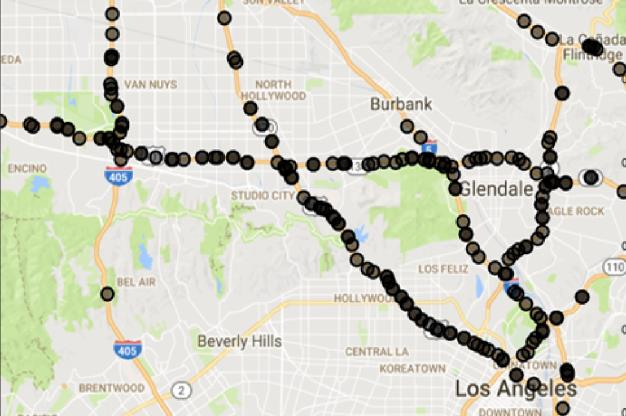}
	\hspace{1.0cm}
	\includegraphics[width=.45\textwidth]{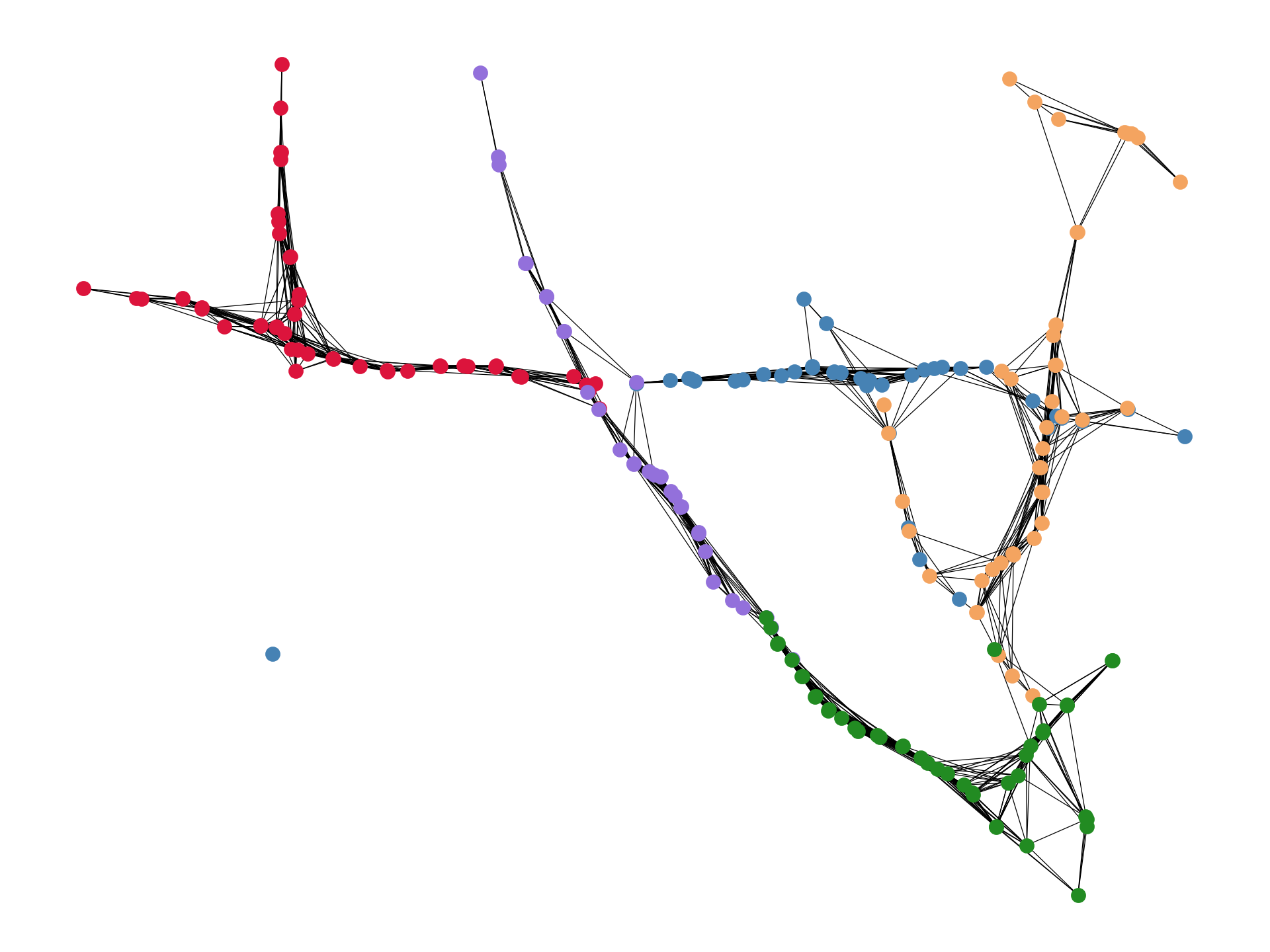}
	\caption{Example of a sample from the METR-LA dataset~\cite{Li2018DiffusionCR}:
		the sensor locations (left) and the corresponding graph partitioned into 5 subgraphs (right).}
	\label{fig:metrla_graph}
\end{figure}

Finally, we evaluate the proposed optimization algorithms on the METR-LA traffic forecasting dataset~\cite{Li2018DiffusionCR}.
This dataset comprises traffic speed measurements from 207 loop detectors deployed on the highway network of Los Angeles County, collected over four months at five-minute intervals.
The learning task is formulated as a node-level spatiotemporal regression problem, in which future traffic conditions are predicted from historical observations on a fixed road network graph.
Each sensor corresponds to a graph node, with traffic measurements represented as node features.
In particular, each node is associated with two input features per time step, i.e., the traffic speed and a time-of-day encoding.
We use the graph provided with the dataset, whose edge weights are obtained by applying a thresholded Gaussian kernel to the pairwise road-network distances between connected sensors (Figure~\ref{fig:metrla_graph}).

We use the preprocessed METR-LA dataset provided by PyTorch Geometric Temporal~\cite{rozemberczki2021pytorch, ockerman2025pgtiscalingspatiotemporalgnns}.
Following the standard forecasting setup of~\cite{Li2018DiffusionCR}, each input sample contains the previous $12$ time steps (corresponding to one hour of observations) and the model predicts the following $12$ time steps.
The traffic measurements are normalized using the training-set statistics.
The resulting $34{,}249$ sliding-window samples are split into $20{,}550$ / $6{,}850$ / $6{,}849$ samples for training, validation, and testing, respectively.

We model the traffic dynamics using the diffusion convolutional recurrent neural network (DCRNN)~\cite{Li2018DiffusionCR}, which combines diffusion convolutions with a GRU-based sequence-to-sequence model.
The employed GNN consists of two DCRNN layers, mapping the two input features to $32$ hidden channels and preserving this width in the second layer, followed by a linear output layer.
Training and evaluation use the masked mean absolute error (MAE) of~\cite{Li2018DiffusionCR}, which ignores missing measurements.

The road-network graph is partitioned using the multilevel Graclus algorithm~\cite{Dhillon2007WeightedGC}, which operates directly on the sparse graph and incorporates edge weights throughout the coarsening process, thereby preserving the road-network connectivity.
Our implementation of the multilevel Graclus algorithm uses the graph matching and coarsening routines provided by PyTorch Geometric~\cite{Fey2019FastGR}.
The initial partitioning of the coarsest graph and the refinement during uncoarsening follow the procedures of~\cite{Dhillon2007WeightedGC}.
Following the heuristic used in~\cite{Dhillon2007WeightedGC}, we coarsen the graph until the coarsest level contains approximately $5P$ nodes, where $P$ denotes the target number of partitions.
The partition labels are then propagated back to the original graph through iterative refinement during uncoarsening.

\subsection{Algorithmic and implementation details}
The GNN models and the proposed training algorithms are implemented in Python~3.12 and rely on PyTorch~2.8~\cite{Paszke2019PyTorchAI}, PyTorch Geometric~\cite{Fey2019FastGR}, and PyTorch Geometric Temporal~\cite{rozemberczki2021pytorch, ockerman2025pgtiscalingspatiotemporalgnns}.
All GNN weights are initialized using Xavier initialization~\cite{Glorot2010UnderstandingDT}, while all biases are initialized to zero.
For batch normalization layers, the scale parameters are initialized to one and the biases to zero.

Unless stated otherwise, both \ddada and \algo use a momentum parameter of $\beta=0.9$.
The coordinate-wise AdaGrad weights are initialized as ${\bm{w}_{0,i}=0.01}$ for all $i=1,\ldots,n$.
The computation of $\gamma_k$ in~\eqref{eq:QP_sub} requires the evaluation of the term $\langle \bm{s}^{S}_{k}, \bm{B}_{k}\bm{s}^{S}_{k}\rangle$.
Rather than explicitly assembling the Hessian $\bm{B}_{k}$, we compute the required Hessian-vector products using PyTorch automatic differentiation.

Graph partitioning is performed once before training, and the resulting partitions are used throughout all DD iterations, since graph partitioning is computationally expensive and the resulting decomposition is deterministic.
In contrast, the coarse graph is regenerated at each epoch via per-subdomain random node subsampling with a prescribed coarsening factor $(c_f)$, defined as the fraction of nodes retained in each subdomain.
This is due to the fact that the coarse graph construction is relatively inexpensive.
Moreover, performing it at every epoch allows for different random coarse representations to be used throughout training, rather than relying on a single fixed subsampling.
The remaining algorithmic parameters, namely the number of outer iterations ($K$), the numbers of global ($K^G$), coarse ($K^C$), and subdomain steps ($K^p$), as well as the number of partitions $P$, are specified separately, as they vary across the experiments.

All experiments were performed on the Jean Zay supercomputer, using GPU nodes equipped with an NVIDIA V100 GPU ($32$~GB of memory) and a 12-core Intel Cascade Lake 6226 CPU.

%% file: sections/05_num_results.tex
\section{Numerical results}
\label{sec:num_results}

In this section, we evaluate the proposed algorithms on the previously described benchmark problems.
To compare the training algorithms under consideration, we report the estimated parallel computational cost $\pazocal{C}$, expressed in units of global-graph optimization steps.
Since the cost of a \sadam optimization step is proportional to the size of the graph on which the MP is performed, a global optimization step on the full graph has a unit cost.
A subdomain optimization step, executed concurrently on all $P$ partitions, has a cost of $\sfrac{1}{P}$, as each partition contains approximately $\sfrac{1}{P}$ of the nodes of the full graph.
Moreover, an optimization step performed on a coarse graph, constructed with the coarsening factor $c_f$, has a cost of $\sfrac{1}{c_f}$.
As a consequence, the estimated computational cost accumulated over $K$ outer iterations is given as
\begin{equation}
	\pazocal{C} := \sum_{k=1}^{K} \qty( K^G + \frac{K^p}{P} + \frac{K^C}{c_f} ),
	\label{eq:comp_cost}
\end{equation}
where $K^G$, $K^p$, and $K^C$ denote the numbers of global, subdomain, and coarse-level optimization steps performed within each outer iteration.

All reported numerical results and plots are averaged over four random seeds.

\subsection{Convergence properties of \sadam}
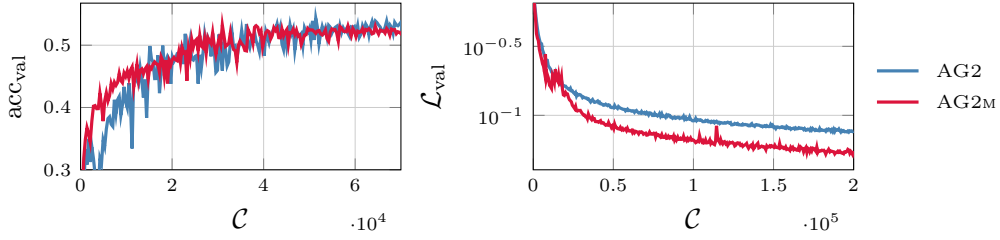
\begin{figure}[tb!]
	\centering
	\input{fig/momentum}
	\caption{Convergence of \sada and \sadam on the CIFAR10 (left) and AirfRANS (right) benchmark problems.
		The results are reported in terms of the validation accuracy (acc$_{\text{val}}$) for CIFAR10 and the validation loss ($\pazocal{L}_{\text{val}}$) for AirfRANS.}
	\label{fig:momentum}
\end{figure}

We begin by investigating the impact of the momentum on the convergence of the single-level algorithm.
To this aim, we compare \sada, i.e., the momentum-free curvature-informed AdaGrad~\cite{Gratton2025NoiseTolerant}, with its momentum-enhanced variant \sadam.
As we can see from Figure~\ref{fig:momentum}, incorporating the momentum improves the convergence.
In particular, on CIFAR10, \sadam converges faster than \sada during the initial phase of the training, while reaching a comparable final validation accuracy.
On AirfRANS, \sadam attains a lower validation loss than \sada for a comparable computational cost.
Based on these results, we employ momentum for all subsequent experiments and considered algorithms.

\subsection{Convergence properties of \ddada}
Next, we investigate the convergence properties of the \ddada algorithm.
In particular, we study the sensitivity of \ddada with respect to the number of global steps $K^G$, the number of subdomain steps $K^p$, and the number of partitions $P$.

\begin{figure}[tb!]
	\centering
	\input{fig/CIFAR10_full_batch}
	\caption{Convergence of \sadam and \ddada with varying $K^G$ for the CIFAR10.
		\ddada is executed with $P=\{2,3,5\}$ (from left to right), and $K^p=|\pazocal{D}|$.}
	\label{fig:CIFAR10_full_batch}
\end{figure}
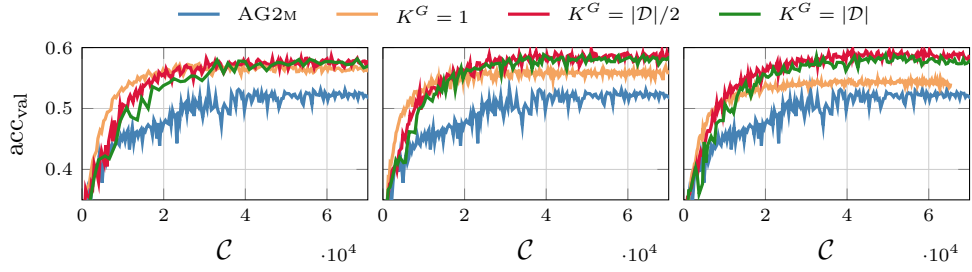

We start by varying the number of global steps $K^G$, while keeping the number of subdomain steps fixed to $K^p=352$, which corresponds to one epoch, i.e., the full pass through the dataset.
Figure~\ref{fig:CIFAR10_full_batch} reports the results obtained for the CIFAR10 benchmark problem with $P \in \{2, 3, 5\}$.
As we can see, using a small number of global steps ($K^G=1$) accelerates the initial phase of the training, since the computational cost of an outer iteration is dominated by the inexpensive subdomain steps.
However, increasing the number of global steps ($K^G \in \{176,352\}$) yields a higher final validation accuracy, with this improvement becoming more pronounced as the number of partitions increases.
Notably, for all considered configurations, \ddada surpasses the \sadam baseline in terms of the validation accuracy reached for a given computational cost $\pazocal{C}$.

\begin{figure}[tb!]
	\centering
	\input{fig/partition_steps}
	\caption{Convergence of \sadam and \ddada ($P=2$, $K^G=|\pazocal{D}|$) with varying $K^p$ for the CIFAR10 (left), for AirfRANS (middle) and for METR-LA (right).}
	\label{fig:partition_steps}
\end{figure}
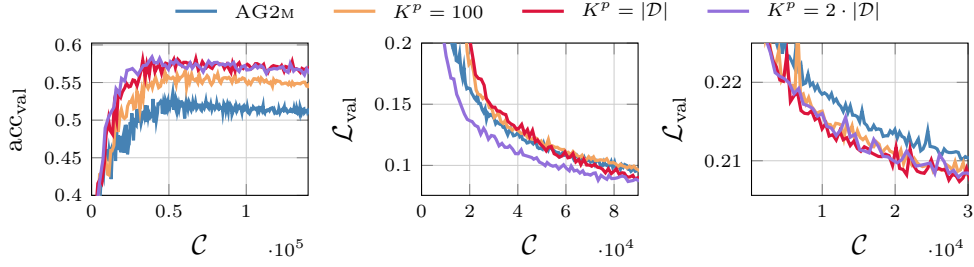

We continue by studying the effect of the number of subdomain steps $K^p$.
Figure~\ref{fig:partition_steps} reports the results obtained for all considered datasets with $P=2$, $K^G=|\pazocal{D}|$, and $K^p \in \{100, |\pazocal{D}|, 2 |\pazocal{D}|\}$.
As we can see, using ${K^p=2 |\pazocal{D}|}$ provides the best or comparable performance across all datasets.
Based on our experimental experience, substantially increasing the number of subdomain steps $K^p$ beyond $2|\pazocal{D}|$ may lead to overfitting of the subgraph models, potentially reducing the benefit of the subsequent global aggregation phase.

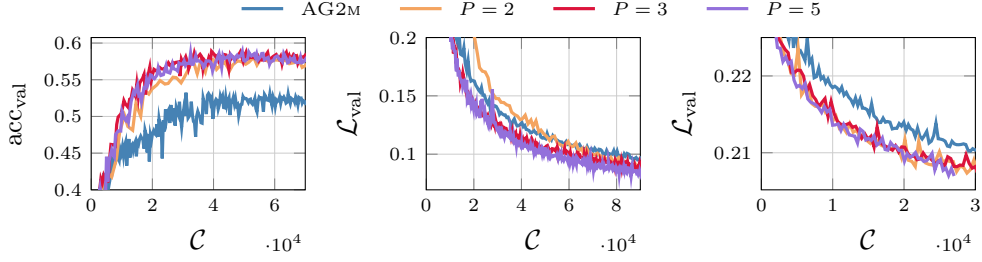
\begin{figure}[tb!]
	\centering
	\input{fig/scaling_DD_wrt_partitions}
	\caption{Convergence of \sadam and \ddada ($K^G=K^p=|\pazocal{D}|$) with increasing number of partitions $P$ for the CIFAR10 (left), AirfRANS (middle), and METR-LA (right).}
	\label{fig:scaling_DD_wrt_partitions}
\end{figure}

Finally, we examine the algorithmic scalability of \ddada with respect to an increasing number of partitions $P$.
As we can see from Figure~\ref{fig:scaling_DD_wrt_partitions}, for all three benchmark problems, \ddada gives rise to a more accurate GNN model than \sadam for a given computational cost $\pazocal{C}$.
We also observe that the convergence of \ddada does not deteriorate with an increasing $P$, which is essential for the parallel scalability of the proposed approach.

\subsection{Convergence properties of \algo}
At the end, we investigate the convergence properties of \algo with respect to the number of coarse steps $K^C$ and the utilized coarsening factor $c_f$.
The performed experiments utilize $K^G=10$.

\begin{figure}[tb!]
	\centering
	\input{fig/CIFAR10_coarse_batch}
	\caption{Convergence of \sadam and \algo ($K^G=10, c_f=2$) with varying $K^C$ for the CIFAR10,  $P=2$ (left), $P=3$ (middle), and $P=5$ (right).}
	\label{fig:CIFAR10_coarse_batch}
\end{figure}
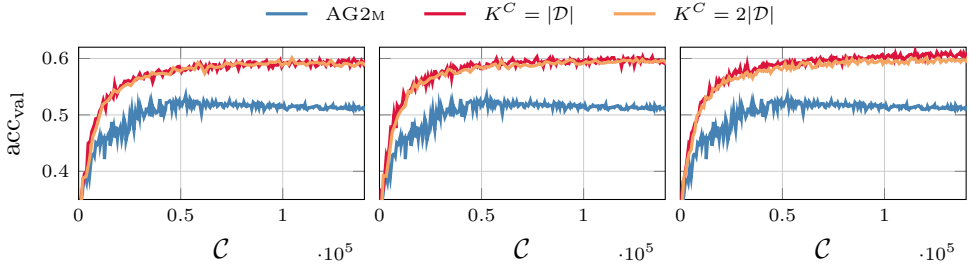

We start by investigating the effect of the number of coarse steps $K^C$ on the convergence behavior of \algo.
The coarse optimization phase replaces a substantial number of expensive full-graph optimization steps with cheaper optimization steps on the coarsened graph.
Figure~\ref{fig:CIFAR10_coarse_batch} reports the results obtained for the CIFAR10 benchmark problem.
As observed, introducing the coarse optimization phase improves the convergence compared to the \sadam baseline, leading to higher final validation accuracy.
Moreover, utilizing $K^C=|\pazocal{D}|$ and $K^C=2 |\pazocal{D}|$ yields comparable performance, which indicates that the number of coarse steps can be reduced without significantly affecting the final solution, thereby improving the parallel scalability.

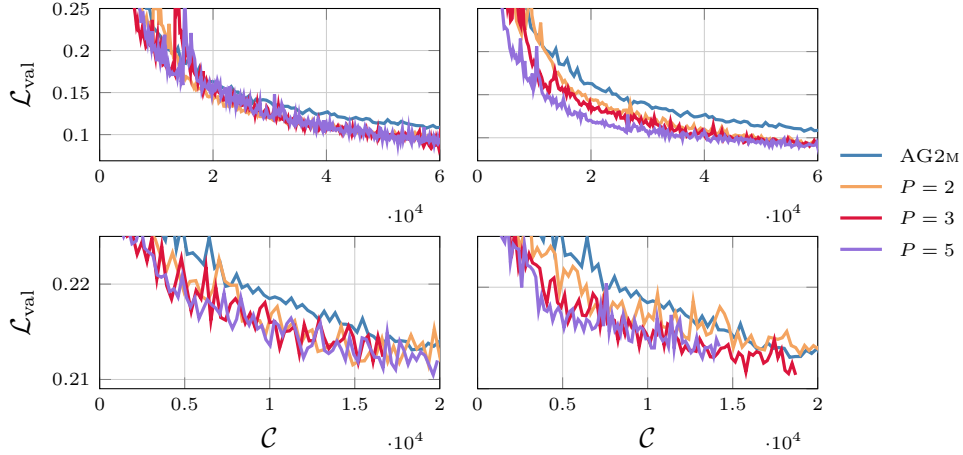
\begin{figure}[tb!]
	\centering
	\input{fig/MA_coarse}
	\caption{Convergence of \sadam and \algo ($K^G=10, K^C=K^p=|\pazocal{D}|$) on the AirfRANS (top) and METR-LA (bottom),  with the coarsening factors $c_f=2$ (left) and $c_f=4$ (right) for different numbers of graph partitions $P$.}
	\label{fig:MA_coarse}
\end{figure}

We next study the behavior of \algo when varying the number of partitions~$P$ and the coarsening factor~$c_f$.
Figure~\ref{fig:MA_coarse} reports the convergence results obtained for the AirfRANS and METR-LA benchmark problems.
For both datasets, increasing the coarsening factor from $c_f=2$ to $c_f=4$ generally improves the convergence, with the difference becoming more significant as the number of partitions increases.
This effect is particularly evident for AirfRANS, where combining a larger coarsening factor with a higher number of partitions leads to the largest performance and accuracy improvements.

\subsection{Comparison between \sadam, \ddada and \algo}
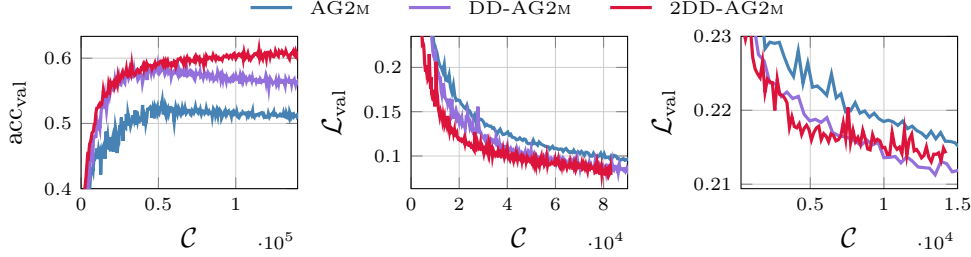
\begin{figure}[tb!]
	\centering
	\input{fig/scaling_TDD}
	\caption{Convergence of \sadam, \ddada ($K^G=K^p=|\pazocal{D}|$, $P=5$) and \algo ($K^G=10, K^C=K^p=|\pazocal{D}|$, $P=5$) on the CIFAR10 (left, $c_f=2$), AirfRANS (middle, $c_f=4$), and METR-LA (right, $c_f=4$).}
	\label{fig:scaling_TDD}
\end{figure}

Finally, we compare the convergence of \sadam, \ddada, and \algo on all three benchmark problems.
Figure~\ref{fig:scaling_TDD} reports the results obtained with $P=5$ partitions.
As we can see, both DD-based methods consistently outperform the \sadam baseline, producing more accurate GNN models for a given computational cost $\pazocal{C}$.
Furthermore, \algo converges faster than \ddada during the early stages of training, and on CIFAR10 it also reaches a notably higher final validation accuracy.
These results highlight the benefits of incorporating the coarse-level steps, which enable the exchange of global information at a fraction of the cost of the full-graph optimization steps.

%% file: fig/momentum.tex
\tikzsetnextfilename{momentum_comparison}
\begin{tikzpicture}

	\pgfplotstableread[col sep=comma]{results/CIFAR10/baseline_validate_acc_plot.csv}\cifarWithout
	\pgfplotstableread[col sep=comma]{results/CIFAR10/baseline_mom_validate_acc_plot.csv}\cifarWith

	\pgfplotstableread[col sep=comma]{results/AirfRANS/baseline_plot.csv}\airWithout
	\pgfplotstableread[col sep=comma]{results/AirfRANS/baseline_no_mom_plot.csv}\airWith

	\begin{groupplot}[
			group style={
					group size=2 by 1,
					horizontal sep=1.75cm,
				},
			width=0.445\linewidth,
			height=0.2895\linewidth, 
			grid=both,
			grid style={line width=.1pt, draw=gray!20},
			major grid style={line width=.2pt, draw=gray!40},
			tick label style={font=\footnotesize},
			label style={font=\footnotesize},
			xlabel near ticks,
			ylabel near ticks,
			title style={font=\footnotesize},
		]


		\nextgroupplot[
			title={},
			xlabel={$\pazocal{C}$},
			ylabel={acc$_{\text{val}}$},
			xmin=0,
			xmax=70000,
			ymin=0.3,
		]

		\addplot[
			very thick,
			color=steelblue,
			mark=none,
		]
		table[x=step,y=Mean]{\cifarWithout};

		\addplot[
			very thick,
			color=crimson,
			mark=none,
		]
		table[x=step,y=Mean]{\cifarWith};


		\nextgroupplot[
			title={},
			xlabel={$\pazocal{C}$},
			ylabel={$\pazocal{L}_{\text{val}}$},
			xmin=0,
			xmax=200000,
			ymode=log,
			ymax=0.65
		]

		\addplot[
			very thick,
			color=steelblue,
			mark=none,
		]
		table[x=step,y=Mean]{\airWithout};

		\addplot[
			very thick,
			color=crimson,
			mark=none,
		]
		table[x=step,y=Mean]{\airWith};

	\end{groupplot}

	\matrix[
		draw=none,
		anchor=west,
		inner sep=0.2em,
		row sep=3pt,
		column sep=4pt,
	] at ([xshift=0.21cm]group c2r1.east) {
		\draw[very thick, color=steelblue] (0,0) -- (0.55,0); &
		\node[font=\footnotesize, anchor=west] {\sada}; \\
		\draw[very thick, color=crimson] (0,0) -- (0.55,0); &
		\node[font=\footnotesize, anchor=west] {\sadam}; \\
	};

\end{tikzpicture}

%% file: fig/CIFAR10_full_batch.tex
\tikzsetnextfilename{CIFAR10_full_batch}

\begin{tikzpicture}

	\pgfplotstableread[col sep=comma]{results/CIFAR10/baseline_mom_validate_acc_plot.csv}\baseline

	\pgfplotstableread[col sep=comma]{results/CIFAR10/P2/fb1_ema_cb0_pre352_plot.csv}\PtwoA
	\pgfplotstableread[col sep=comma]{results/CIFAR10/P2/fb176_ema_cb0_pre352_plot.csv}\PtwoB
	\pgfplotstableread[col sep=comma]{results/CIFAR10/P2/fb352_ema_cb0_pre352_plot.csv}\PtwoC

	\pgfplotstableread[col sep=comma]{results/CIFAR10/P3/fb1_ema_cb0_pre352_plot.csv}\PthreeA
	\pgfplotstableread[col sep=comma]{results/CIFAR10/P3/fb176_ema_cb0_pre352_plot.csv}\PthreeB
	\pgfplotstableread[col sep=comma]{results/CIFAR10/P3/fb352_ema_cb0_pre352_plot.csv}\PthreeC

	\pgfplotstableread[col sep=comma]{results/CIFAR10/P5/fb1_ema_cb0_pre352_plot.csv}\PfiveA
	\pgfplotstableread[col sep=comma]{results/CIFAR10/P5/fb176_ema_cb0_pre352_plot.csv}\PfiveB
	\pgfplotstableread[col sep=comma]{results/CIFAR10/P5/fb352_ema_cb0_pre352_plot.csv}\PfiveC

	\begin{groupplot}[
			group style={
					group size=3 by 1,
					horizontal sep=0.2cm,
				},
			width=0.41\linewidth,
			height=0.275\linewidth,
			grid=both,
			ymin=0.35,
			grid style={line width=.1pt, draw=gray!20},
			major grid style={line width=.2pt, draw=gray!40},
			tick label style={font=\footnotesize},
			label style={font=\footnotesize},
			xlabel near ticks,
			ylabel near ticks,
			title style={font=\footnotesize},
			xmin=0,
			ymax=0.6,
			xlabel={$\pazocal{C}$},
			xmax=70000,
		]

		\nextgroupplot[
			title={},
			ylabel={acc$_{\text{val}}$},
		]

		\addplot[very thick, color=steelblue, mark=none] table[x=step,y=Mean]{\baseline};
		\addplot[very thick, color=sandybrown, mark=none] table[x=step,y=Mean]{\PtwoA};
		\addplot[very thick, color=crimson, mark=none] table[x=step,y=Mean]{\PtwoB};
		\addplot[very thick, color=forestgreen, mark=none] table[x=step,y=Mean]{\PtwoC};

		\nextgroupplot[
			ylabel={},
			title={},
			yticklabel=\empty,
		]

		\addplot[very thick, color=steelblue, mark=none] table[x=step,y=Mean]{\baseline};
		\addplot[very thick, color=sandybrown, mark=none] table[x=step,y=Mean]{\PthreeA};
		\addplot[very thick, color=crimson, mark=none] table[x=step,y=Mean]{\PthreeB};
		\addplot[very thick, color=forestgreen, mark=none] table[x=step,y=Mean]{\PthreeC};

		\nextgroupplot[
			ylabel={},
			title={},
			yticklabel=\empty,
		]

		\addplot[very thick, color=steelblue, mark=none] table[x=step,y=Mean]{\baseline};
		\addplot[very thick, color=sandybrown, mark=none] table[x=step,y=Mean]{\PfiveA};
		\addplot[very thick, color=crimson, mark=none] table[x=step,y=Mean]{\PfiveB};
		\addplot[very thick, color=forestgreen, mark=none] table[x=step,y=Mean]{\PfiveC};

	\end{groupplot}

	\matrix[
		draw=none,
		anchor=south,
		inner sep=0.2em,
		column sep=4pt,
	] at ([yshift=0.15cm]group c2r1.north) {
		\draw[very thick, color=steelblue] (0,0) -- (0.55,0);   &
		\node[font=\footnotesize, anchor=west] {\sadam};        & [8pt]
		\draw[very thick, color=sandybrown] (0,0) -- (0.55,0);  &
		\node[font=\footnotesize, anchor=west] {$K^G=1$};       & [8pt]
		\draw[very thick, color=crimson] (0,0) -- (0.55,0);     &
		\node[font=\footnotesize, anchor=west] {$K^G=|\pazocal{D}|/2$};     & [8pt]
		\draw[very thick, color=forestgreen] (0,0) -- (0.55,0); &
		\node[font=\footnotesize, anchor=west] {$K^G=|\pazocal{D}|$};             \\
	};

\end{tikzpicture}

%% file: fig/partition_steps.tex
\tikzsetnextfilename{partition_steps}

\begin{tikzpicture}

	\pgfplotstableread[col sep=comma]{results/CIFAR10/baseline_mom_validate_acc_plot.csv}\cifarBaseline
	\pgfplotstableread[col sep=comma]{results/CIFAR10/P2/fb352_ema_cb0_pre100_plot.csv}\cifarPreOneHundred
	\pgfplotstableread[col sep=comma]{results/CIFAR10/P2/fb352_ema_cb0_pre352_plot.csv}\cifarPreDataset
	\pgfplotstableread[col sep=comma]{results/CIFAR10/P2/fb352_ema_cb0_pre704_plot.csv}\cifarPreTwoDataset

	\pgfplotstableread[col sep=comma]{results/AirfRANS/baseline_plot.csv}\airfBaseline
	\pgfplotstableread[col sep=comma]{results/AirfRANS/P2/fb200_cb0_pre100_plot.csv}\airfPreOneHundred
	\pgfplotstableread[col sep=comma]{results/AirfRANS/P2/fb200_cb0_pre200_plot.csv}\airfPreDataset
	\pgfplotstableread[col sep=comma]{results/AirfRANS/P2/fb200_cb0_pre400_plot.csv}\airfPreTwoDataset

	\pgfplotstableread[col sep=comma]{results/METR_LA/baseline_plot.csv}\metrBaseline
	\pgfplotstableread[col sep=comma]{results/METR_LA/P2/fb322_cb0_pre100_plot.csv}\metrPreOneHundred
	\pgfplotstableread[col sep=comma]{results/METR_LA/P2/fb322_cb0_pre322_plot.csv}\metrPreDataset
	\pgfplotstableread[col sep=comma]{results/METR_LA/P2/fb322_cb0_pre644_plot.csv}\metrPreTwoDataset

	\begin{groupplot}[
			group style={
					group size=3 by 1,
					horizontal sep=1.5cm,
				},
			width=0.34\linewidth,
			height=0.275\linewidth,
			grid=both,
			grid style={line width=.1pt, draw=gray!20},
			major grid style={line width=.2pt, draw=gray!40},
			tick label style={font=\footnotesize},
			label style={font=\footnotesize},
			xlabel near ticks,
			ylabel near ticks,
			title style={font=\footnotesize},
			xlabel={$\pazocal{C}$},
		]


		\nextgroupplot[
			title={},
			ylabel={acc$_{\text{val}}$},
			xmin=0,
			xmax=140000,
			ymin=0.40,
		]

		\addplot[very thick, color=steelblue, mark=none] table[x=step,y=Mean]{\cifarBaseline};
		\addplot[very thick, color=sandybrown, mark=none] table[x=step,y=Mean]{\cifarPreOneHundred};
		\addplot[very thick, color=crimson, mark=none] table[x=step,y=Mean]{\cifarPreDataset};
		\addplot[very thick, color=mediumpurple, mark=none] table[x=step,y=Mean]{\cifarPreTwoDataset};


		\nextgroupplot[
			title={},
			ylabel={$\pazocal{L}_{\text{val}}$},
			xmin=0,
			ymax=0.2,
			xmax=89975,
		]

		\addplot[very thick, color=steelblue, mark=none] table[x=step,y=Mean]{\airfBaseline};
		\addplot[very thick, color=sandybrown, mark=none] table[x=step,y=Mean]{\airfPreOneHundred};
		\addplot[very thick, color=crimson, mark=none] table[x=step,y=Mean]{\airfPreDataset};
		\addplot[very thick, color=mediumpurple, mark=none] table[x=step,y=Mean]{\airfPreTwoDataset};


		\nextgroupplot[
			title={},
			ylabel={$\pazocal{L}_{\text{val}}$},
			xmin=300,
			ymax=0.225,
			ytick={0.21,0.22,0.23},
			xmax=30000,
		]

		\addplot[very thick, color=steelblue, mark=none] table[x=step,y=Mean]{\metrBaseline};
		\addplot[very thick, color=sandybrown, mark=none] table[x=step,y=Mean]{\metrPreOneHundred};
		\addplot[very thick, color=crimson, mark=none] table[x=step,y=Mean]{\metrPreDataset};
		\addplot[very thick, color=mediumpurple, mark=none] table[x=step,y=Mean]{\metrPreTwoDataset};

	\end{groupplot}

	\matrix[
		draw=none,
		anchor=south,
		inner sep=0.2em,
		column sep=4pt,
	] at ([yshift=0.15cm]group c2r1.north) {
		\draw[very thick, color=steelblue] (0,0) -- (0.55,0);                    &
		\node[font=\footnotesize, anchor=west] {\sadam};                         & [8pt]
		\draw[very thick, color=sandybrown] (0,0) -- (0.55,0);                   &
		\node[font=\footnotesize, anchor=west] {$K^p=100$};                      & [8pt]
		\draw[very thick, color=crimson] (0,0) -- (0.55,0);                      &
		\node[font=\footnotesize, anchor=west] {$K^p=|\pazocal{D}|$}; & [8pt]
		\draw[very thick, color=mediumpurple] (0,0) -- (0.55,0);                 &
		\node[font=\footnotesize, anchor=west] {$K^p=2\cdot|\pazocal{D}|$};   \\
	};

\end{tikzpicture}

%% file: fig/scaling_DD_wrt_partitions.tex
\tikzsetnextfilename{scaling_DD_wrt_partitions}
\begin{tikzpicture}

	\pgfplotstableread[col sep=comma]{results/CIFAR10/baseline_mom_validate_acc_plot.csv}\cifarBaseline
	\pgfplotstableread[col sep=comma]{results/CIFAR10/P2/fb352_ema_cb0_pre352_plot.csv}\cifarPtwo
	\pgfplotstableread[col sep=comma]{results/CIFAR10/P3/fb352_ema_cb0_pre352_plot.csv}\cifarPthree
	\pgfplotstableread[col sep=comma]{results/CIFAR10/P5/fb352_ema_cb0_pre352_plot.csv}\cifarPfive

	\pgfplotstableread[col sep=comma]{results/AirfRANS/baseline_plot.csv}\airBaseline
	\pgfplotstableread[col sep=comma]{results/AirfRANS/P2/fb200_cb0_pre200_plot.csv}\airPtwo
	\pgfplotstableread[col sep=comma]{results/AirfRANS/P3/fb200_cb0_pre200_plot.csv}\airPthree
	\pgfplotstableread[col sep=comma]{results/AirfRANS/P5/fb200_cb0_pre200_plot.csv}\airPfive
	\pgfplotstableread[col sep=comma]{results/AirfRANS/P8/fb200_cb0_pre200_plot.csv}\airPeight

	\pgfplotstableread[col sep=comma]{results/METR_LA/baseline_plot.csv}\metrBaseline
	\pgfplotstableread[col sep=comma]{results/METR_LA/P2/fb322_cb0_pre322_plot.csv}\metrPtwo
	\pgfplotstableread[col sep=comma]{results/METR_LA/P3/fb322_cb0_pre322_plot.csv}\metrPthree
	\pgfplotstableread[col sep=comma]{results/METR_LA/P5/fb322_cb0_pre322_plot.csv}\metrPfive

	\begin{groupplot}[
			group style={
					group size=3 by 1,
					horizontal sep=1.6cm,
				},
			width=0.3375\linewidth,
			height=0.275\linewidth,
			grid=both,
			grid style={line width=.1pt, draw=gray!20},
			major grid style={line width=.2pt, draw=gray!40},
			tick label style={font=\footnotesize},
			label style={font=\footnotesize},
			xlabel near ticks,
			ylabel near ticks,
			title style={font=\footnotesize},
			xlabel={$\pazocal{C}$},
		]


		\nextgroupplot[
			title={},
			ylabel={acc$_{\text{val}}$},
			xmin=0,
			ymin=0.4,
			xmax=70000,
		]

		\addplot[very thick, color=steelblue, mark=none] table[x=step,y=Mean]{\cifarBaseline};
		\addplot[very thick, color=sandybrown, mark=none] table[x=step,y=Mean]{\cifarPtwo};
		\addplot[very thick, color=crimson, mark=none] table[x=step,y=Mean]{\cifarPthree};
		\addplot[very thick, color=mediumpurple, mark=none] table[x=step,y=Mean]{\cifarPfive};


		\nextgroupplot[
			title={},
			ylabel={$\pazocal{L}_{\text{val}}$},
			xmin=0,
			ymax=0.20,
			xmax=89975,
		]

		\addplot[very thick, color=steelblue, mark=none] table[x=step,y=Mean]{\airBaseline};
		\addplot[very thick, color=sandybrown, mark=none] table[x=step,y=Mean]{\airPtwo};
		\addplot[very thick, color=crimson, mark=none] table[x=step,y=Mean]{\airPthree};
		\addplot[very thick, color=mediumpurple, mark=none] table[x=step,y=Mean]{\airPfive};


		\nextgroupplot[
			title={},
			ylabel={$\pazocal{L}_{\text{val}}$},
			xmin=0,
			xmax=30000,
			ymax=0.225,
			ytick={0.21,0.22,0.23},
		]

		\addplot[very thick, color=steelblue, mark=none] table[x=step,y=Mean]{\metrBaseline};
		\addplot[very thick, color=sandybrown, mark=none] table[x=step,y=Mean]{\metrPtwo};
		\addplot[very thick, color=crimson, mark=none] table[x=step,y=Mean]{\metrPthree};
		\addplot[very thick, color=mediumpurple, mark=none] table[x=step,y=Mean]{\metrPfive};

	\end{groupplot}

	\matrix[
		draw=none,
		anchor=south,
		inner sep=0.2em,
		column sep=4pt,
	] at ([yshift=0.15cm]group c2r1.north) {
		\draw[very thick, color=steelblue] (0,0) -- (0.55,0);    &
		\node[font=\footnotesize, anchor=west] {\sadam};         & [8pt]
		\draw[very thick, color=sandybrown] (0,0) -- (0.55,0);   &
		\node[font=\footnotesize, anchor=west] {$P=2$};          & [8pt]
		\draw[very thick, color=crimson] (0,0) -- (0.55,0);      &
		\node[font=\footnotesize, anchor=west] {$P=3$};          & [8pt]
		\draw[very thick, color=mediumpurple] (0,0) -- (0.55,0); &
		\node[font=\footnotesize, anchor=west] {$P=5$};                  \\
	};

\end{tikzpicture}

%% file: fig/CIFAR10_coarse_batch.tex
\tikzsetnextfilename{CIFAR10_coarse_batch}

\begin{tikzpicture}

	\pgfplotstableread[col sep=comma]{results/CIFAR10/baseline_mom_validate_acc_plot.csv}\baseline

	\pgfplotstableread[col sep=comma]{results/CIFAR10/P2/fb10_ema_cb352_cl1_pre352_ab10_plot.csv}\PtwoB
	\pgfplotstableread[col sep=comma]{results/CIFAR10/P2/fb10_ema_cb704_cl1_pre352_ab10_plot.csv}\PtwoC

	\pgfplotstableread[col sep=comma]{results/CIFAR10/P3/fb10_ema_cb352_cl1_pre352_ab10_plot.csv}\PthreeB
	\pgfplotstableread[col sep=comma]{results/CIFAR10/P3/fb10_ema_cb704_cl1_pre352_ab10_plot.csv}\PthreeC

	\pgfplotstableread[col sep=comma]{results/CIFAR10/P5/fb10_ema_cb352_cl1_pre352_ab10_plot.csv}\PfiveB
	\pgfplotstableread[col sep=comma]{results/CIFAR10/P5/fb10_ema_cb704_cl1_pre352_ab10_plot.csv}\PfiveC

	\begin{groupplot}[
			group style={
					group size=3 by 1,
					horizontal sep=0.2cm,
				},
			width=0.41\linewidth,
			height=0.275\linewidth,
			grid=both,
			ymin=0.35,
			grid style={line width=.1pt, draw=gray!20},
			major grid style={line width=.2pt, draw=gray!40},
			tick label style={font=\footnotesize},
			label style={font=\footnotesize},
			xlabel near ticks,
			ylabel near ticks,
			title style={font=\footnotesize},
			xmin=0,
			ymax=0.62,
			xlabel={$\pazocal{C}$},
			xmax=140000,
		]

		\nextgroupplot[
			title={},
			ylabel={acc$_{\text{val}}$},
		]

		\addplot[very thick, color=steelblue, mark=none] table[x=step,y=Mean]{\baseline};
		\addplot[very thick, color=crimson, mark=none] table[x=step,y=Mean]{\PtwoB};
		\addplot[very thick, color=sandybrown, mark=none] table[x=step,y=Mean]{\PtwoC};

		\nextgroupplot[
			ylabel={},
			title={},
			yticklabel=\empty,
		]

		\addplot[very thick, color=steelblue, mark=none] table[x=step,y=Mean]{\baseline};
		\addplot[very thick, color=crimson, mark=none] table[x=step,y=Mean]{\PthreeB};
		\addplot[very thick, color=sandybrown, mark=none] table[x=step,y=Mean]{\PthreeC};

		\nextgroupplot[
			ylabel={},
			title={},
			yticklabel=\empty,
		]

		\addplot[very thick, color=steelblue, mark=none] table[x=step,y=Mean]{\baseline};
		\addplot[very thick, color=crimson, mark=none] table[x=step,y=Mean]{\PfiveB};
		\addplot[very thick, color=sandybrown, mark=none] table[x=step,y=Mean]{\PfiveC};

	\end{groupplot}

	\matrix[
		draw=none,
		anchor=south,
		inner sep=0.2em,
		column sep=4pt,
	] at ([yshift=0.15cm]group c2r1.north) {
		\draw[very thick, color=steelblue] (0,0) -- (0.55,0);  &
		\node[font=\footnotesize, anchor=west] {\sadam};       & [8pt]
		\draw[very thick, color=crimson] (0,0) -- (0.55,0);    &
		\node[font=\footnotesize, anchor=west] {$K^C=|\pazocal{D}|$};    & [8pt]
		\draw[very thick, color=sandybrown] (0,0) -- (0.55,0); &
		\node[font=\footnotesize, anchor=west] {$K^C=2|\pazocal{D}|$};            \\
	};

\end{tikzpicture}

%% file: fig/MA_coarse.tex
\tikzsetnextfilename{MA_coarse}
\begin{tikzpicture}

	\pgfplotstableread[col sep=comma]{results/AirfRANS/baseline_plot.csv}\airBaseline

	\pgfplotstableread[col sep=comma]{results/AirfRANS/P2/fb10_cb200_cl1_pre200_ab10_plot.csv}\airPtwoLone
	\pgfplotstableread[col sep=comma]{results/AirfRANS/P3/fb10_cb200_cl1_pre200_ab10_plot.csv}\airPthreeLone
	\pgfplotstableread[col sep=comma]{results/AirfRANS/P5/fb10_cb200_cl1_pre200_ab10_plot.csv}\airPfiveLone

	\pgfplotstableread[col sep=comma]{results/AirfRANS/P2/fb10_cb200_cl2_pre200_ab10_plot.csv}\airPtwoLtwo
	\pgfplotstableread[col sep=comma]{results/AirfRANS/P3/fb10_cb200_cl2_pre200_ab10_plot.csv}\airPthreeLtwo
	\pgfplotstableread[col sep=comma]{results/AirfRANS/P5/fb10_cb200_cl2_pre200_ab10_plot.csv}\airPfiveLtwo

	\pgfplotstableread[col sep=comma]{results/METR_LA/baseline_plot.csv}\metrBaseline

	\pgfplotstableread[col sep=comma]{results/METR_LA/P2/fb10_cb322_cl1_pre322_ab10_plot.csv}\metrPtwoLone
	\pgfplotstableread[col sep=comma]{results/METR_LA/P3/fb10_cb322_cl1_pre322_ab10_plot.csv}\metrPthreeLone
	\pgfplotstableread[col sep=comma]{results/METR_LA/P5/fb10_cb322_cl1_pre322_ab10_plot.csv}\metrPfiveLone

	\pgfplotstableread[col sep=comma]{results/METR_LA/P2/fb10_cb322_cl2_pre322_ab10_plot.csv}\metrPtwoLtwo
	\pgfplotstableread[col sep=comma]{results/METR_LA/P3/fb10_cb322_cl2_pre322_ab10_plot.csv}\metrPthreeLtwo
	\pgfplotstableread[col sep=comma]{results/METR_LA/P5/fb10_cb322_cl2_pre322_ab10_plot.csv}\metrPfiveLtwo

	\begin{groupplot}[
			group style={
					group size=2 by 2,
					horizontal sep=0.5cm,
					vertical sep=1.cm,
				},
			width=0.465\linewidth,
			height=0.275\linewidth,
			grid=both,
			grid style={line width=.1pt, draw=gray!20},
			major grid style={line width=.2pt, draw=gray!40},
			tick label style={font=\footnotesize},
			label style={font=\footnotesize},
			xlabel near ticks,
			ylabel near ticks,
			title style={font=\footnotesize},
			xmin=0,
			ymax=0.25,
		]


		\nextgroupplot[
			title={},
			xlabel={},
			ylabel={$\pazocal{L}_{\text{val}}$},
			xmax=60000,
		]

		\addplot[very thick, color=steelblue, mark=none] table[x=step,y=Mean]{\airBaseline};
		\addplot[very thick, color=sandybrown, mark=none] table[x=step,y=Mean]{\airPtwoLone};
		\addplot[very thick, color=crimson, mark=none] table[x=step,y=Mean]{\airPthreeLone};
		\addplot[very thick, color=mediumpurple, mark=none] table[x=step,y=Mean]{\airPfiveLone};


		\nextgroupplot[
			title={},
			xlabel={},
			ylabel={},
			yticklabel=\empty,
			xmax=60000,
		]

		\addplot[very thick, color=steelblue, mark=none] table[x=step,y=Mean]{\airBaseline};
		\addplot[very thick, color=sandybrown, mark=none] table[x=step,y=Mean]{\airPtwoLtwo};
		\addplot[very thick, color=crimson, mark=none] table[x=step,y=Mean]{\airPthreeLtwo};
		\addplot[very thick, color=mediumpurple, mark=none] table[x=step,y=Mean]{\airPfiveLtwo};


		\nextgroupplot[
			title={},
			xlabel={$\pazocal{C}$},
			ylabel={$\pazocal{L}_{\text{val}}$},
			xmax=20000,
			ymax=0.225,
			ytick={0.21,0.22,0.23},
		]

		\addplot[very thick, color=steelblue, mark=none] table[x=step,y=Mean]{\metrBaseline};
		\addplot[very thick, color=sandybrown, mark=none] table[x=step,y=Mean]{\metrPtwoLone};
		\addplot[very thick, color=crimson, mark=none] table[x=step,y=Mean]{\metrPthreeLone};
		\addplot[very thick, color=mediumpurple, mark=none] table[x=step,y=Mean]{\metrPfiveLone};


		\nextgroupplot[
			title={},
			xlabel={$\pazocal{C}$},
			ylabel={},
			yticklabel=\empty,
			xmax=20000,
			ymax=0.225,
			ytick={0.21,0.22,0.23},
		]

		\addplot[very thick, color=steelblue, mark=none] table[x=step,y=Mean]{\metrBaseline};
		\addplot[very thick, color=sandybrown, mark=none] table[x=step,y=Mean]{\metrPtwoLtwo};
		\addplot[very thick, color=crimson, mark=none] table[x=step,y=Mean]{\metrPthreeLtwo};
		\addplot[very thick, color=mediumpurple, mark=none] table[x=step,y=Mean]{\metrPfiveLtwo};

	\end{groupplot}

	\matrix[
		draw=none,
		anchor=west,
		inner sep=0.2em,
		row sep=3pt,
		column sep=4pt,
	] at ([xshift=0.2cm, yshift=-0.55cm]group c2r1.south east) {
		\draw[very thick, color=steelblue] (0,0) -- (0.55,0);    &
		\node[font=\footnotesize, anchor=west] {\sadam};           \\
		\draw[very thick, color=sandybrown] (0,0) -- (0.55,0);   &
		\node[font=\footnotesize, anchor=west] {$P=2$};            \\
		\draw[very thick, color=crimson] (0,0) -- (0.55,0);      &
		\node[font=\footnotesize, anchor=west] {$P=3$};            \\
		\draw[very thick, color=mediumpurple] (0,0) -- (0.55,0); &
		\node[font=\footnotesize, anchor=west] {$P=5$};            \\
	};

\end{tikzpicture}

%% file: fig/scaling_TDD.tex
\tikzsetnextfilename{scaling_TDD}
\begin{tikzpicture}

	\pgfplotstableread[col sep=comma]{results/CIFAR10/baseline_mom_validate_acc_plot.csv}\cifarBaseline
	\pgfplotstableread[col sep=comma]{results/CIFAR10/P5/fb352_ema_cb0_pre352_plot.csv}\cifarDD
	\pgfplotstableread[col sep=comma]{results/CIFAR10/P5/fb10_ema_cb352_cl1_pre352_ab10_plot.csv}\cifarTDD

	\pgfplotstableread[col sep=comma]{results/AirfRANS/baseline_plot.csv}\airfBaseline
	\pgfplotstableread[col sep=comma]{results/AirfRANS/P5/fb200_cb0_pre200_plot.csv}\airfDD
	\pgfplotstableread[col sep=comma]{results/AirfRANS/P5/fb10_cb200_cl2_pre200_ab10_plot.csv}\airfTDD

	\pgfplotstableread[col sep=comma]{results/METR_LA/baseline_plot.csv}\metrBaseline
	\pgfplotstableread[col sep=comma]{results/METR_LA/P5/fb322_cb0_pre322_plot.csv}\metrDD
	\pgfplotstableread[col sep=comma]{results/METR_LA/P5/fb10_cb322_cl2_pre322_ab10_plot.csv}\metrTDD

	\begin{groupplot}[
			group style={
					group size=3 by 1,
					horizontal sep=1.5cm,
				},
			width=0.34\linewidth,
			height=0.275\linewidth,
			grid=both,
			grid style={line width=.1pt, draw=gray!20},
			major grid style={line width=.2pt, draw=gray!40},
			tick label style={font=\footnotesize},
			label style={font=\footnotesize},
			xlabel near ticks,
			ylabel near ticks,
			title style={font=\footnotesize},
			xlabel={$\pazocal{C}$},
		]


		\nextgroupplot[
			title={},
			ylabel={acc$_{\text{val}}$},
			xmin=0,
			xmax=140000,
			ymin=0.40,
		]

		\addplot[very thick, color=steelblue, mark=none] table[x=step,y=Mean]{\cifarBaseline};
		\addplot[very thick, color=mediumpurple, mark=none] table[x=step,y=Mean]{\cifarDD};
		\addplot[very thick, color=crimson, mark=none] table[x=step,y=Mean]{\cifarTDD};


		\nextgroupplot[
			title={},
			ylabel={$\pazocal{L}_{\text{val}}$},
			xmin=0,
			ymax=0.235,
			xmax=89975,
			yticklabel style={
					/pgf/number format/fixed,
					/pgf/number format/precision=2
				},
		]

		\addplot[very thick, color=steelblue, mark=none] table[x=step,y=Mean]{\airfBaseline};
		\addplot[very thick, color=mediumpurple, mark=none] table[x=step,y=Mean]{\airfDD};
		\addplot[very thick, color=crimson, mark=none] table[x=step,y=Mean]{\airfTDD};


		\nextgroupplot[
			title={},
			ylabel={$\pazocal{L}_{\text{val}}$},
			xmin=300,
			xmax=15000,
			ymax=0.23,
			ytick={0.21,0.22,0.23},
		]

		\addplot[very thick, color=steelblue, mark=none] table[x=step,y=Mean]{\metrBaseline};
		\addplot[very thick, color=mediumpurple, mark=none] table[x=step,y=Mean]{\metrDD};
		\addplot[very thick, color=crimson, mark=none] table[x=step,y=Mean]{\metrTDD};

	\end{groupplot}

	\matrix[
		draw=none,
		anchor=south,
		inner sep=0.2em,
		column sep=4pt,
	] at ([yshift=0.15cm]group c2r1.north) {
		\draw[very thick, color=steelblue] (0,0) -- (0.55,0);    &
		\node[font=\footnotesize, anchor=west] {\sadam};         & [8pt]
		\draw[very thick, color=mediumpurple] (0,0) -- (0.55,0); &
		\node[font=\footnotesize, anchor=west] {\ddada};         & [8pt]
		\draw[very thick, color=crimson] (0,0) -- (0.55,0);      &
		\node[font=\footnotesize, anchor=west] {\algo};                  \\
	};

\end{tikzpicture}

%% file: sections/06_conclusion.tex
\section{Conclusion}
\label{sec:conclusion}
In this work, we proposed \ddada and \algo, two DD-based algorithms for the scalable training of GNNs.
Both algorithms are built upon \sadam, a curvature-informed AdaGrad method with momentum~\cite{Gratton2025NoiseTolerant}.
They combine optimization steps on graph partitions with optimization steps performed either on the full graph (\ddada) or on a coarse graph obtained by random subsampling within each subdomain (\algo).
To evaluate the proposed methods, we conducted numerical experiments on graph classification, node-level regression, and spatiotemporal forecasting tasks using different GNN architectures and graph partitioning strategies.
Our results demonstrate that \ddada and \algo require up to $4$–$8\times$ fewer optimization steps than \sadam to achieve the same predictive performance.
At the same time, \ddada and \algo improve the predictive performance of GNNs by up to $22\%$ for the same number of optimization steps.
Moreover, the performance remained stable as the number of graph partitions increased, demonstrating the algorithmic scalability of the proposed algorithms.

We foresee several extensions of the presented work.
In particular, we plan to develop a distributed-memory implementation to translate the observed reduction in computational cost into a wall-clock speedup.
We also aim to investigate the use of overlapping subdomains, alternative coarsening strategies, and multilevel extensions with more than two levels.
Finally, we plan to design asynchronous variants, including their practical implementation and analysis of their convergence properties.